\documentclass[oneside, 11pt]{article}
\usepackage{setspace}
\usepackage{titlesec}
\usepackage[utf8]{inputenc}
\usepackage{tabularx}	
\usepackage{etoolbox}
\usepackage{mathtools}
\usepackage{amsmath}
\usepackage{amssymb} 
\usepackage{amsfonts}
\usepackage{mathrsfs} 
\usepackage{amsthm,pifont}
\usepackage{graphicx}
\usepackage{bbm}
\usepackage[colorlinks=true, pdfstartview=FitV, linkcolor=blue,
citecolor=blue, urlcolor=blue, pagebackref=false]{hyperref}
\usepackage[nottoc]{tocbibind}
\usepackage{xcolor}
\usepackage{tikz}
\usepackage{setspace}
\usepackage{cite}
\usepackage{braket}
\usepackage{leftidx,tensor}
\usepackage{dsfont}
\usepackage{calc}
\usepackage{caption}
\usepackage{subcaption}
\usepackage{color}
\usetikzlibrary{decorations.markings}
\usetikzlibrary{shapes.geometric}
\usetikzlibrary{patterns}
\tikzstyle{vertex}=[circle, draw, inner sep=0pt, minimum size=6pt]

\usepackage[utf8]{inputenc}
\usepackage[english]{babel}
\usepackage{enumitem}

\allowdisplaybreaks

\newtheorem{theorem}{Theorem}[section]
\makeatletter
\patchcmd{\ttlh@hang}{\parindent\z@}{\parindent\z@\leavevmode}{}{}
\patchcmd{\ttlh@hang}{\noindent}{}{}{}
\makeatother
\titleformat*{\section}{\large\bfseries}
\titleformat*{\subsection}{\small\bfseries}
\titleformat*{\subsubsection}{\small\bfseries}
\titleformat*{\paragraph}{\small\bfseries}
\titleformat*{\subparagraph}{\small\bfseries}

\newcommand{\N}{\mathbb{N}}
\newcommand{\R}{\mathbb{R}}
\newcommand{\Z}{\mathbb{Z}}
\newcommand{\E}{\mathbb{E}}

\newcommand{\p}{\mathbb{P}}
\newcommand{\md}{\ensuremath{\mathrm{d}}}
\newcommand{\eps}{\varepsilon}

\newcommand{\cD}{\mathcal{D}}
\newcommand{\cE}{\mathcal{E}}
\newcommand{\cF}{\mathcal{F}}

\newcommand{\cP}{\mathcal{P}}

\newcommand{\cX}{\mathcal{X}}

\newtheorem{lemma}[theorem]{Lemma}
\newtheorem{remark}[theorem]{Remark}
\newtheorem{proposition}[theorem]{Proposition}

\newtheorem{corollary}[theorem]{Corollary}

\usepackage{geometry}
\begin{document}

	\title{On self-similar random quotients of the line}
	
	\author{Johannes B\"aumler\thanks{Mathematical Institute, University of Koblenz; \url{jbaeumler@uni-koblenz.de}}
    \and
    Tom Hutchcroft\thanks{Department of Mathematics, Princeton University  and The Division of Physics, Mathematics and Astronomy, California Institute of Technology; \url{t.hutchcroft@princeton.edu}}
    }

	\maketitle
	
	\begin{center}
		\parbox{13cm}{ \textbf{Abstract.} 
			A surprisingly wide variety of important random fractal metric spaces in modern probability theory can be viewed as self-similar random quotients of the real line, including the continuum random tree, Liouville quantum gravity surfaces, the trace of Brownian motion in $\R^d$, and self-similar ($s=2d$) long-range percolation on $\R$. We prove that the metric and resistance growth dimensions are well-defined in a certain sense for every random metric space in a class containing all of the above examples. Our results are new for the trace of Brownian motion.
		}
	\end{center}
	
	
	
	\setstretch{1.1}

	\section{Introduction}\label{sec:introduction}

	\subsection{Renormalization of distances in self-similar random quotients of the real line}

    In this paper we study the metric and spectral properties of \emph{self-similar random quotients of the real line}. Such a random quotient is defined topologically by its set of identification points.
    Let $\Omega$ denote the set of unordered pairs of distinct elements of $\R$. We are primarily interested in random closed\footnote{Restricting to closed sets is done only to give us a nice measure space to work with and is not important for any of our results.} subsets $\cX$ of $\Omega$ 
    satisfying the following axioms:
    \begin{enumerate}[label=(A\arabic*),ref=A\arabic*]
		\item[(I)] \textbf{Independence of non-overlapping intervals: }\label{assu:indep} For each $a_1 \leq b_1 \leq a_2 \leq b_2 \leq \cdots \leq a_n \leq b_n$, the sets $(\{\{x,y\} \in \cX :a_i \leq x,y \leq b_i\})_{i=1}^n$ are independent.        
		\item[(T)] \textbf{Translation invariance:}
        \label{assu:translation} The random set $a+\cX := \left\{\{x+a,y+a\}: \{x,y\} \in \cX\right\}$ has the same law as $\cX$ for each $a\in \R$.
        \item[(S)] \label{assu:self-sim} \textbf{Scale invariance:} The random set $\lambda \cX=\{\{\lambda x,\lambda y\} : \{x,y\} \in \cX\}$ has the same distribution as $\cX$ for every $\lambda>0$.
\end{enumerate}
Note that we do \emph{not} require $\cX$ to be an equivalence relation, and passing to the transitive closure of $\cX$ may break the independence assumption (I).

\begin{figure}[t!]
\centering
\includegraphics[clip, trim = 0cm 0.9cm 0cm 0cm, width=\textwidth]{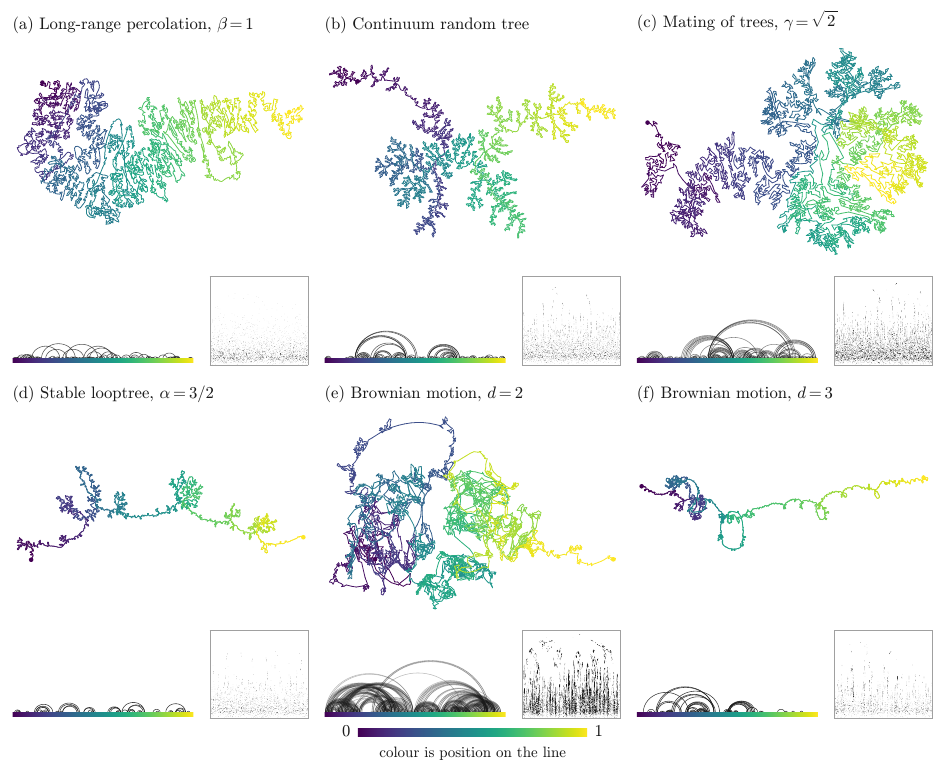}
\caption{Simulations of several of the models falling under the hypotheses of our main theorems. In each case, we have plotted a spring-electrical embedding of the relevant graph metric into the plane (top; note that only the $\Z$ edges are shown, so that each space is depicted with a distinguished space-filling curve), the pair process $\cX$ shown as an arc diagram (bottom left), and the pair process $\cX$ with each point $\{x,y\}\in \cX$ plotted as $((x+y)/2N, \log (N|x-y|))$ (bottom right). Each simulation was carried out on an interval of length $N=2^{12}=4096$.}
\label{fig:simulations_pretty}
\end{figure}

\begin{figure}
\centering
\includegraphics[clip, trim = 0cm 0.8cm 0cm 0cm, width=\textwidth]{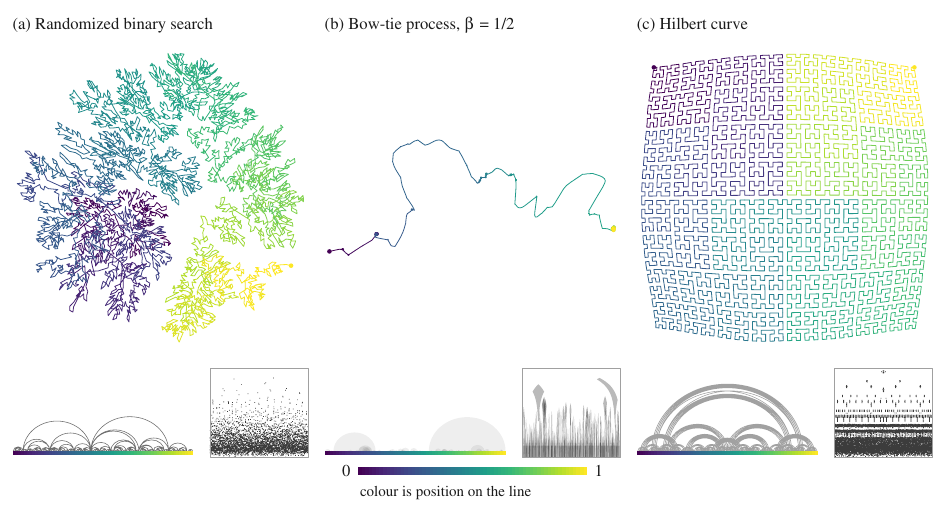}
\caption{Left and center: The counterexamples from Section~\ref{subsec:counterexamples_to_stronger_statements} plotted in the same manner as in Figure~\ref{fig:simulations_pretty}. Right: The Euclidean plane can be thought of as a quotient of the real line by a fractal set of identifications defined by a space filling curve (here the Hilbert curve); this set $\cX$ does not satisfy our axioms (T) or (S).}
\label{fig:counterexamples}
\end{figure}

A surprisingly large array of important random fractal metric spaces\footnote{No prior knowledge of any of these examples is needed to read the paper.} in modern probability can be interpreted as quotients $\R/\cX$ or $[0,1]/\cX$ for random sets satisfying the axioms above. Note that each of these examples is endowed with a canonical choice of \emph{space filling curve} defined from its quotient structure.
\begin{enumerate}
\item If we take $\cX$ to be a Poisson point process on $\Omega$ with intensity $\beta |x-y|^{-2}$ for some $\beta>0$, the associated random quotients $\R/\cX$  and $[0,1]/\cX$ can be thought of in terms of \emph{long-range percolation} at the self-similar point $s=2d$. Relevant references include  \cite{benjamini2001diameter,coppersmith2002diameter,baumler2023distances,baumler2023polynomial,ding2013distances,ding2023uniqueness,ding2025uniqueness} for the metric structure and \cite{MR1896880,ding2025polynomial,baumler2023recurrence,ding2025growth} for the behaviour of random walk.
\item Let $(B_t)_{t\in \R}$ be standard  Brownian motion in $\R^d$ and let $\cX = \{\{s,t\} : B_t=B_s\}$ be the associated set of self-intersections. The quotients $\R/\cX$ and $[0,1]/\cX$ encode the trace of Brownian motion, the metric properties of which are poorly understood in two and three dimensions \cite{burdzy1990rigorous}; related non-rigorous and numerical studies include \cite{zhou2026fractals,manna1989diffusion,dekeyser1987diffusion}. 
(In dimensions $d\geq 4$ Brownian motion is almost surely simple, but one can study analogous questions for the discrete random walk, which exhibit logarithmic corrections to scaling in four dimensions \cite{croydon2009random,croydon2023scaling,shiraishi2012exact,shiraishi2026graph}.)
 More generally one can take $\cX = \{\{s,t\}:L_t=L_s$ or $\lim_{u\uparrow s}L_{u}\}$ for any (strictly) $\alpha$-stable L\'evy process $(L_t)_{t\in \R}$.
\item Let $(B_t)_{t\in \R}$ be Brownian motion in $\R$ and let $\cX = \{\{s,t\} : B_t=B_s$ and $B_r \geq B_t$ for every $r$ between $s$ and $t\}$. The quotient $[0,1]/\cX$ is closely related to the \emph{continuum random tree} (CRT) of fixed mass, which is obtained by the same construction with $B$ replaced by a Brownian excursion, while $\R/\cX$ is the infinite version of the CRT \cite{CRT1,aldous1991continuum,aldous1993continuum}. Important variations on this construction include \emph{stable trees} \cite{duquesne2002random,duquesne2005hausdorff} and \emph{stable looptrees} \cite{curien2014random,archer2021brownian}, which are constructed using certain stable L\'evy processes instead of Brownian motion and describe the scaling limits of critical random trees with heavy-tailed degree distributions.
\item Let $(B^1_t,B^2_t)$ be a pair of possibly correlated Brownian motions in $\R$ (so that $(B^1_t,B^2_t)$ is a possibly anisotropic Brownian motion in $\R^2$), and let 
$\cX = \{\{s,t\} : B_t^1=B_s^1$ and $B_r^1 \geq B_t^1$ for every $r$ between $s$ and $t\}\cup \{\{s,t\} : B_t^2=B_s^2$ and $B_r^2 \geq B_t^2$ for every $r$ between $s$ and $t\}$. 
The mating of trees theory of Duplantier, Miller, and Sheffield interprets the quotient $\R/\cX$ in terms of \emph{Liouville quantum gravity} (LQG) surfaces \cite{duplantier2021liouville} (see also \cite{gwynne2019mating}), with different universality classes of LQG surfaces encoded by Brownian motions with different correlation constants. 
Further references on the geometric and spectral properties of LQG surfaces and related discrete models include \cite{gwynne2019mating,gwynne2020mating,gwynne2019harmonic,gwynne2021random,gwynne2020anomalous,miller2021axiomatic,ding2023uniquenessLQG,ding2020fractal,gwynne2019distance,lee2023relations}.
Important variations on these constructions include \emph{matings of stable looptrees} \cite[Section 1.3]{duplantier2021liouville}, which describe the scaling limits of critical random clusters in statistical mechanics models on random maps and in which one again replaces the Brownian motion in the construction of LQG surfaces with appropriate stable L\'evy processes.
\end{enumerate}
See Figure~\ref{fig:simulations_pretty} for simulations, and Figure~\ref{fig:counterexamples}(c) for an example of a process $\cX$ that does \emph{not} satisfy these axioms. A more detailed literature review is given in Section~\ref{subsec:exponent_values_and_related_literature}.

In each case above, trying to define a metric on $\R/\cX$ or $[0,1]/\cX$ naively via
\[
\inf \Biggl\{ \sum_{i=0}^n |u_{i}-v_i| : n \geq 0, u_0 = x, v_n =y, \{v_i,u_{i+1}\} \in \cX \text{ for every $0\leq i \leq n-1$}\Biggr\}
\]
makes all distances equal to zero almost surely\footnote{Indeed this holds whenever $\cX \neq \emptyset$ as a consequence of our main theorems.}, so one must instead use some renormalization procedure to define the metric, 
for example by giving each pair in $\cX$ a small positive cost $\eps$ to cross and rescaling the resulting metric by an appropriately chosen function of $\eps$.  
Formally, given a random set of pairs $\cX$ as above, we define for each interval $I \subseteq \R$ and $\eps>0$ the metric
\begin{multline*}
d_I^{(\eps)}(x,y) = 
\inf \Biggl\{ \eps n+ \sum_{i=0}^n |u_{i}-v_i| : n \geq 0, u_0 = x, v_{n} =y, u_i,v_i\in I\\ \text{ for every $0\leq i \leq n$, and }  \{v_i,u_{i+1}\} \in \cX \text{ for every $0\leq i \leq n-1$}\Biggr\}.
\end{multline*}
In order to define a limit of these metrics, it is important to understand their dependence on $\eps$. For many examples of interest we expect the system to have a well-defined \emph{dimension} in the sense that distances decay like a well-defined power of $\eps$. The existence of such a dimension has been proven via different methods for many of the examples above, with relevant references including 
\cite{baumler2023distances,baumler2023polynomial,ding2013distances,gwynne2020mating,gwynne2019distance},
but remained open for e.g.\ the trace of Brownian motion in two and three dimensions. 

\medskip

In this paper we give a soft, general argument guaranteeing the existence of a well-defined distance growth exponent for all processes $\cX$ satisfying the above conditions. We also prove a similar theorem for the  \emph{effective resistance}, which we state in Section~\ref{subsec:effective_resistances}.

\begin{theorem}
\label{thm:main}
Let $\mathcal{X}$ be a random closed subset of $\Omega=\left\{ \{x,y\}\subset \R : x \neq y \right\}$. If $\cX$ satisfies \emph{(I)}, \emph{(T)}, and \emph{(S)} then there exists $\theta\in [0,1]$  such that 
\[
\E \left[ d^{(\eps)}_{[0,1]}(0,1) \right] = \eps^{1-\theta+o(1)} \qquad\text{ and } \qquad
\E \left[ \sup_{x,y \in [0,1]} d^{(\eps)}_{[0,1]}(x,y) \right] = \eps^{1-\theta+o(1)}
\]
as $\eps\downarrow 0$. Moreover, if $\mathcal{X} \neq \varnothing$ with positive probability then $\theta<1$.
\end{theorem}

The Brownian case of our theorems is closely related to a conjecture of Burdzy and Lawler \cite{burdzy1990rigorous} (which concerns discrete random walk rather than Brownian motion).

The proof of the theorem also shows that the relevant distances are of order $\eps^{1-\theta+o(1)}$ with good probability and are exponentially unlikely to be larger than a quantity of order $\eps^{1-\theta+o(1)}$ by a large multiplicative factor (though this quantity may have a  different $o(1)$ term in its exponent than the mean does). The proof does not show that the relevant distances have order at least $\eps^{1-\theta+o(1)}$ with high probability, and indeed this is not true at the level of generality of the theorem as we discuss in Section~\ref{subsec:counterexamples_to_stronger_statements}.

\medskip

\noindent \textbf{Discrete large-scale formulation.} 
Note that Theorem~\ref{thm:main} can also be formulated as a result about the metric $d^{(1)}$ at large scales, with the distributional equality
\begin{equation}
\label{eq:basic_scaling_identity}
\Bigl(d^{(\eps)}_{[0,1]}(x,y):x,y\in [0,1]\Bigr) \overset{d}{=} \Bigl(\eps d^{(1)}_{[0,1/\eps]}(x/\eps,y/\eps):x,y\in [0,1]\Bigr)
\end{equation}
holding by scale invariance of $\cX$.
We now state an alternative, more discrete large-scale formulation of our results. Let $\cX$ be a random subset of $\Omega$ as above. We use the process $\cX$ to define a \emph{graph} with vertex set $\Z$ by setting two integers $x$ and $y$ to be adjacent if either $|x-y|=1$ or $|x-y|>1$ and there exist points $u\in [x,x+1]$ and $v\in [y,y+1]$ such that $\{u,v\}\in \cX$. We denote the resulting graph metric on the discrete interval $I\subseteq \Z$ by $D_I$, denoting the graph itself by $G_I$. We call these graphs the \textbf{block graphs} associated to $\cX$. We also write $D_{[a,b]}=D_{\{a,\ldots,b-1\}}$ and $G_{[a,b]} = G_{\{a,\ldots,b-1\}}$ for integers $a<b$, noting that this graph is determined by the restriction of $\cX$ to $[a,b]$.

\begin{theorem}
\label{thm:main_discrete}
Let $\mathcal{X}$ be a random closed subset of $\Omega=\left\{ \{x,y\}\subset \R : x \neq y \right\}$. If $\cX$ satisfies  \emph{(I)}, \emph{(T)}, and \emph{(S)}  then there exists $\theta\in [0,1]$  such that 
\[
\E \left[ D_{[0,n]}(0,n-1) \right] = n^{\theta+o(1)} \qquad\text{ and } \qquad
\E \left[ \sup_{x,y \in \{0,\ldots,n-1\}} D_{[0,n]}(x,y) \right] = n^{\theta+o(1)}
\]
as $n\to\infty$. This exponent $\theta$ coincides with that of Theorem~\ref{thm:main} and hence satisfies $\theta<1$ when $\cX$ is non-empty with positive probability.
\end{theorem}

As before, the proof of this theorem also establishes that the relevant distances have order at least $n^{\theta-o(1)}$ with good probability but not necessarily with high probability.

\begin{remark}
\emph{
Assuming one can exactly compute the distribution of the block graph $G_{[0,n]}$ for every $n\geq 1$, the proof of Theorem~\ref{thm:main_discrete} leads naturally to an algorithm that outputs certificates for rigorous upper bounds on $\theta$ that converge to the true value as $n\to \infty$. Unfortunately the running time of this algorithm appears to make it completely impractical to use.
}
\end{remark}

\begin{remark}
\emph{As we were in the final stages of preparing this work, the preprint \cite{shiraishigraphdistance2026} appeared, in which Shiraishi gives an independent proof that the chemical distance exponent is well-defined for the trace of 3d Brownian motion. Shiraishi's paper establishes the stronger result that distances in the block graph
$G_{[0,n]}$ are of order $n^{\theta+o(1)}$ with high probability, and 
establishes that $\theta$ is equal to half the minimal Hausdorff dimension of a curve connecting $B(0)$ and $B(1)$ in $B([0,1])$. On the other hand, \cite{shiraishigraphdistance2026} does not study the effective resistance or take the unifying approach considered here. Both our paper and \cite{shiraishigraphdistance2026} leave open the problem of extending their results to the trace of simple random walk.
}
\end{remark}

\subsection{Effective resistances}
\label{subsec:effective_resistances}

We now discuss extensions of our main theorems to \emph{effective resistances}. 
Here we recall that the effective resistance between two vertices $u$ and $v$ in a finite graph $G$ is defined by
\[
	\mathscr{R}_\mathrm{eff} (u\leftrightarrow v; G)^{-1} = \deg(u) \mathbf{P}_u(\text{hit $v$ before returning to $u$}),
\]
where $\mathbf{P}_u$ denotes the law of the simple random walk on $G$ started at $u$. Equivalent expressions for the effective resistance of a more metric-geometric flavour are recalled in Section~\ref{sec:effective_resistances}. Understanding the effective resistance is usually the first step towards understanding the behaviour of random walk on (random) fractals; see \cite{KumagaiBook} for an introduction and \cite{lee2023relations} for the most general results in the setting of random graphs.

We formulate our effective resistance results only for the block graph, with the large number of parallel edges causing problems with the other formulation.
Given the process $\cX \subseteq \Omega$ and the associated block graphs $G_I$ as defined above, we define $R_I(x,y)$ to be the effective resistance between $x$ and $y$ in the block graph $G_I$.
 As before, we also write $R_{[0,n]}$ to denote this quantity in the case $I=\{0,\ldots,n-1\}$. Our results on the effective resistance will require the following mild additional hypothesis:
\begin{itemize}
	\item[(L)] For each $\eps>0$, the maximal degree of $G_{[0,n]}$ is bounded by $n^{\eps}$ with probability tending to $1$ as $n\to \infty$. 
\end{itemize}
 It is easily proven that the trace of Brownian motion satisfies (L) in two and three dimensions as we prove in Section~\ref{subsec:L_brownian}.
The assumption (L) does \emph{not} hold for the trace of the $\alpha$-stable L\'evy process in one dimension with $1<\alpha\leq 2$ or for the $\alpha$-stable tree with $\alpha<2$.

\begin{theorem}
\label{thm:main_discrete_other_metrics}
Let $\mathcal{X}$ be a random closed subset of $\Omega=\left\{ \{x,y\}\subset \R : x \neq y \right\}$. If $\cX$ satisfies  \emph{(I)}, \emph{(T)}, \emph{(S)}, and \emph{(L)} then there exists $\theta_\mathrm{res}\in [0,1]$  such that 
\[
\E \left[ R_{[0,n]}(0,n-1) \right] = n^{\theta_\mathrm{res}+o(1)} \qquad\text{ and } \qquad
\E \left[ \sup_{x,y \in \{0,\ldots,n-1\}} R_{[0,n]}(x,y) \right] = n^{\theta_\mathrm{res}+o(1)}
\]
as $n\to\infty$.
\end{theorem}

The proof of this theorem shows without invoking hypothesis (L) that the growth exponent is always well-defined for the \emph{vertex extremal length}, which is equivalent to the effective resistance up to a subpolynomial factor when (L) holds.
Note that effective resistances are always smaller than graph distances, so that $\theta_\mathrm{res}\leq \theta<1$ when $\cX$ is non-trivial in the sense that $\p(\cX \neq \emptyset)>0$. As with Theorem~\ref{thm:main} and \ref{thm:main_discrete}, the Brownian case of this theorem is very closely related to conjectures of Burdzy and Lawler \cite{burdzy1990rigorous}.


\subsection{Counterexamples to stronger statements}
\label{subsec:counterexamples_to_stronger_statements}

We now discuss counterexamples to some natural strengthenings of our main theorems.

\medskip

\noindent \textbf{Failure of tightness.} Given the statements of Theorems~\ref{thm:main} and \ref{thm:main_discrete}, it is natural to wonder whether the relevant distances are the correct order \emph{with high probability}. The following example shows that this does not follow from assumptions (I), (T), and (S). Consider a Poisson process $\cX$ on $\Omega$ with intensity $\beta |x-y|^{-2}$, and let $\cX^*$ be the set of pairs $s<t$ such that there exists a pair $u<v$ in $\cX$ such that $u\leq s \leq (v+u)/2 \leq t \leq v$. We call this the \textbf{bow-tie process}. It is easily verified that this process $\cX^*$ satisfies the three axioms (I), (T), and (S). Moreover, if $\beta<1$ then the expected number of \emph{cut points} of $G_{[0,n]}$, i.e.\ vertices $i$ for which there is no edge $\{a,b\}$ with $a<i<b$, grows polynomially in $n$ for both $\cX$ and $\cX^*$ (indeed, $i$ is a cut point for $\cX$ with probability $((i+1)(n-i)/n)^{-\beta}$ and for $\cX^*$ with probability $2^{-\beta}((i+1)(n-i))^{-\beta/2}$, both of which sum to order $n^{1-\beta}$; a second moment calculation shows that the number of cut points has this order this good probability also), which implies that the chemical distance exponent $\theta$ associated to $\cX^*$ must be positive. On the other hand, with positive probability there exists an arc of $\cX$ covering the entire interval $[0,n]$ whose midpoint lies in $[0,n]$, in which case the diameter of the block graph $G_{[0,n]}$ associated to $\cX^*$ is at most $2$. This example can be made less degenerate by filling the interior of each arc with a bipartite variant of the binary search example below instead of the complete bipartite graph.

\medskip
\noindent \textbf{Small distances.} The set of all pairs of distinct points in $\R$ satisfies (I), (T), and (S) and is easily seen to satisfy the conclusions of Theorems~\ref{thm:main} and \ref{thm:main_discrete} with $\theta=0$, with the associated block graphs $G_I$ being complete graphs.
We now show that it is possible for $\theta$ to be equal to zero in much less degenerate examples than this, with the associated graphs not only being sparse but having exponential tails on their degrees. (Note that $\theta_\mathrm{res}=0$ is rather more common, and holds for e.g.\ the mating of trees example \cite{gwynne2021random}.)
Let $\cD$ be a Poisson process on $\R \times (0,\infty)$ whose intensity has density $g(x,y)=y^{-2}$ with respect to Lebesgue, and define $\cX$ to be the set of pairs $u<v$ in $\R$ for which there exist $a,b \in (0,\infty)$ such that $(u,a)$ and $(v,b)$ both belong to $\cD$ and there does not exist a point $(w,c)$ of $\cD$ with $w$ between $u$ and $v$ and $c > \min\{a,b\}$. It is easily verified that $\cX$ satisfies (I), (T), and (S). We claim the associated metric has $\theta=0$ and that the associated graph $G_\Z$ has finite degrees, with the probability that the origin has degree at least $k$ being exponentially small in $k$. We will see that this random graph encodes a random binary search of the interval.

We begin by analyzing the metric. We will prove that the diameter of $G_{[0,n]}$ is $O(\log n)$ with high probability.
Fix $m \in \{0,1,\ldots,n-1\}$.
 Let $a_0=0$, $b_0=n$, and let $(x_0,y_0)$ be the point of $\cD$ with the largest second coordinate among those points with first coordinate in $[0,n]$, this point being well-defined since the density $g(x,y)=y^{-2}$ gives finite measure to points with second coordinate at least $h$ for each $h>0$. The point $x_0$ is uniformly distributed on $[0,n]$ by the Poisson splitting property. 
If $x_0>m$ we set $a_1=0$ and $b_1=x_0$, while if $x_0<m$ we set $a_1=x_0$ and $b_1=n$.
 Recursively, for each $i\geq 1$ let $(x_{i},y_{i})$ be the point of $\cD$ with the largest second coordinate among those points with first coordinate in $(a_i,b_i)$ and define $(a_{i+1},b_{i+1})$ by $(a_{i+1},b_{i+1})=(a_i,x_i)$ if $x_i>m$ and 
 $(a_{i+1},b_{i+1})=(x_i,b_i)$ if $x_i<m$. It is easily proven via induction and the Poisson splitting property that
$m\in [a_i,b_i]$ for every $i\geq 0$ and that
  $x_{i+1}$ is uniform on $[a_{i+1},b_{i+1}]$ conditional on $x_0,\ldots,x_i$ for each $i\geq 0$. The definition of $\cX$ ensures that $\lfloor x_{i+1}\rfloor$ is equal or adjacent to $\lfloor x_i \rfloor$ in the graph $G_{[0,n]}$ for every $i\geq 0$.
Meanwhile, the random variables $U_i=(x_i-a_i)/(b_i-a_i)$ are i.i.d.\ uniform on $[0,1]$ with $|a_{i+1}-b_{i+1}|/|a_{i}-b_{i}| \leq \max\{U_i,1-U_i\}$, so Hoeffding's inequality implies that there exists a constant $C$ such that if $k = \lceil C \log (n+1)\rceil$ then
\begin{equation*}
\p(|a_{k}-b_{k}| < 1) \geq	\p\left( \sum_{i=0}^{k-1} \log \max\{U_i,1-U_i\} < - \log (n+1)  \right) 
\geq 1 - (n+1)^{-2}.
\end{equation*}
When this event holds, we must have that $\lfloor x_{k} \rfloor \in \{m-1,m\}$ and hence that $m$ has graph distance at most $k+1$ from $\lfloor x_0 \rfloor$. Since $x_0$ did not depend on the choice of $m$, it follows from this and a union bound that the graph $G_{[0,n]}$ has diameter at most $2k+2 = 2\lceil C \log (n+1)\rceil + 2$ with high probability and hence that $\theta=0$ as claimed.

We next control the degree of the origin in the graph $G_\Z$. For each $k\in \Z$ let $Y_k=\max\{y : (x,y) \in \cD,\, x \in [k,k+1)\}$, so that the different values $Y_k$ are independent and identically distributed. Observe that for an integer $k\geq 1$ to be adjacent to $0$ it must satisfy $Y_k = \max\{Y_i : 0<i\leq k\}$. If we define the sequence of record times $(R_i)_{i=1}^\infty$ recursively by $R_1=1$ and $R_{i+1}=\min\{k>R_i : Y_k > \max\{Y_j : 1\leq j < k\}\}$, it is a property of arbitrary i.i.d.\ sequences that $\p(R_i \leq e^{c_1 i}) \leq Ce^{-c_2 i}$ for every $i\geq 1$ and some universal positive constants $c_1,c_2,C>0$. (Indeed, this follows from R\'enyi's observation that if the measure is non-atomic then the random variables $\mathbf{1}(Y_i = \max_{1\leq j\leq i} Y_j)$ are \emph{independent} with $\p(Y_i = \max_{1\leq j\leq i} Y_j)=1/j$ \cite{renyi1962theorie}, together with standard concentration estimates for sums of independent Bernoulli random variables.) 
%
%
On the other hand, in order for the origin to be connected to any vertex $m\in \{1,2,\ldots\}$, there must not exist any points $(x_1,y_1)\in \cD$ with $y_1 \geq Y_0$ and $x_1\in [1,m)$. Given $Y_0$, the probability that this occurs is 
\[
	\exp\left[- (m-1)\int_{Y_0}^\infty y^{-2} d y\right] = \exp\left[- \frac{(m-1)}{Y_0}\right].
\]
Letting $\deg(0)^+$ be the number of vertices $m\geq 1$ such that the origin is adjacent to $m$ in $G_\Z$, it follows easily that there exist positive constants $c_3$ and $C'$ such that
\begin{multline*}
	\p(\deg(0)^+ \geq k) \leq \p(R_k \leq e^{c_1 k}) + \p(Y_0 \geq e^{c_1 k/2}) \\+ \p(Y_0 \leq e^{c_1k/2}, \text{$0$ is adjacent to some $m\geq e^{c_1k}$})
	\leq C'e^{-c_3 k}
\end{multline*}
for every $k\geq 1$. It follows by symmetry and a union bound that $\p(\deg(0) \geq k)$ is exponentially small in $k$ also.

\subsection{Exponent values and related literature}
\label{subsec:exponent_values_and_related_literature}

\begin{figure}[t!]
\centering
\includegraphics[width=\textwidth]{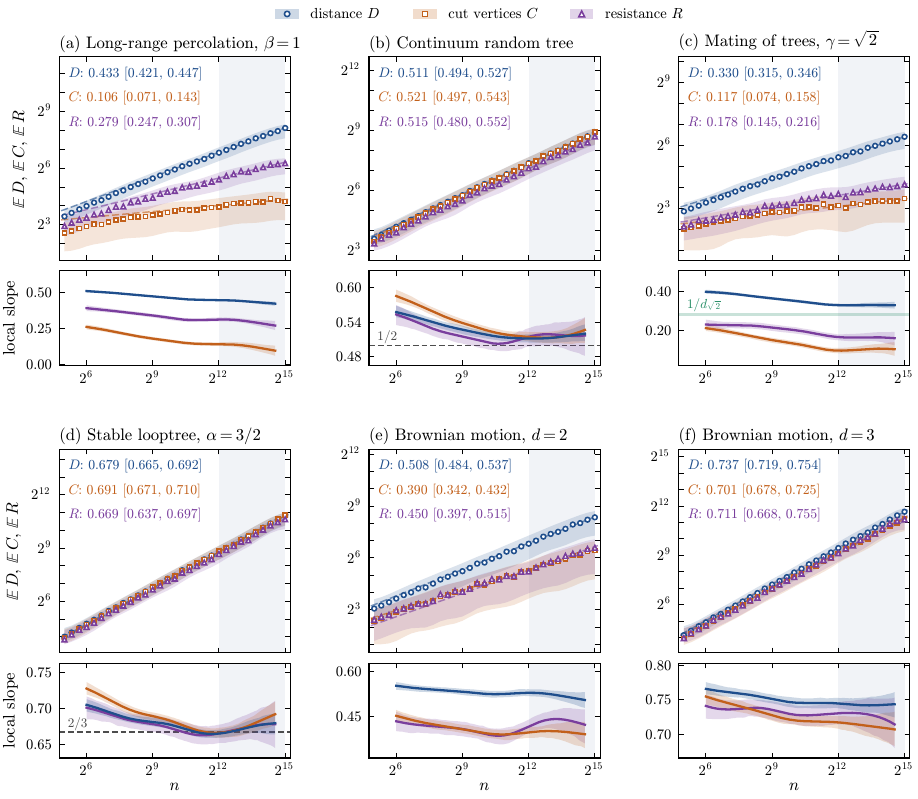}
\caption{Numerical estimates of the growth dimension of the chemical distance, effective resistance, and number of cutpoints in the same examples considered in Figure~\ref{fig:simulations_pretty}. All exponent estimates are computed using only samples of size between $2^{12}$ and $2^{15}$. The confidence intervals reported here are for the empirical values at the finite system sizes studied and do not necessarily contain the true limiting exponent values (as is the case for (a) where the cut point exponent is zero and (c) where the exponents for the resistances and cut points are $0$ \cite{gwynne2021random}). Simulations (b)-(d) use constructions based on discrete random walks rather than the exactly scale-invariant constructions studied in this paper.}
\label{fig:numerical_estimates}
\end{figure}

We now briefly survey what is known regarding the values of the exponents whose existence is proven here. In each case, we can bound distances from below by the number of cut points of $G_{[0,n]}$ (as defined in Section~\ref{subsec:counterexamples_to_stronger_statements}), with the relevant exponent often being easier to understand than the corresponding chemical distance and effective resistance exponents. The same lower bound holds for effective resistances up to a factor of the maximal degree, since consecutive cut points are joined in series and the effective resistance between two distinct vertices is at least the reciprocal of the degree of either one. Heuristically, $\theta$ should be the reciprocal of the dimension of the limiting metric space, which is known in several of our examples, but we do not justify this rigorously for every example we consider. We also report the outcome of simulations estimating the exponents in Figure~\ref{fig:numerical_estimates}.

\medskip

 \noindent \textbf{Brownian motion.}
For $d=1$ it is not hard to see that $\theta=1/2$, with the random metric space $d^{(1)}_{[0,t]}$ being a close approximation to the interval $[\min_{0\leq s \leq t} B(s), \max_{0\leq s \leq t} B(s)]$. 
More generally, it is straightforward to prove that $\theta\geq 1/2$ in every dimension by comparison with the Euclidean displacement. In dimensions $d=2,3$ we can also bound $\theta$ from below by the growth exponent $1-\xi_d(1,1)/2$ of the number of Brownian cut times, where $\xi_d(1,1)$ is the Brownian intersection exponent \cite{lawler1996hausdorff,burdzy1990non}. It should also be possible to bound $\theta$ from above by $\beta_d/2$, where $\beta_d$ is the growth exponent of loop-erased random walk (but this requires proof\footnote{The bound $\theta \leq \beta/2$ is proven for the trace of 3d Brownian motion in the parallel independent work of Shiraishi \cite{shiraishigraphdistance2026}.} since we have defined $\theta$ in terms of discretized Brownian motion instead of simple random walk). In two dimensions $\xi_2(1,1)=5/4$ \cite{lawler2001values2}, so that the resulting lower bound of $3/8$ is weaker than the displacement bound, while $\beta_2=5/4$ \cite{Ken00}. In three dimensions $1/2 <\xi_3(1,1) \leq 0.9999$ and $\beta_3\in[3/2,5/3]$ \cite{lawler1999loop,lawler1989intersections,shiraishi2018growth,cai2026quadratic,liu2026explicit}. These bounds should yield that $1/2\leq \theta\leq 5/8$ in two dimensions and $0.50005\leq \theta \leq 5/6$ in three dimensions. For the effective resistance, the cut-time bound continues to apply by Lemma~\ref{lem:brownian_L} and the discussion above, so that $\theta_\mathrm{res}\geq 1-\xi_d(1,1)/2$, yielding that $\theta_\mathrm{res}\geq 3/8$ in two dimensions and $\theta_\mathrm{res}\geq 0.50005$ in three dimensions. It is conjectured that $\theta=1/2$ in two dimensions, in agreement with the numerical studies \cite{dekeyser1987diffusion,manna1989diffusion,zhou2026fractals} and with Figure~\ref{fig:numerical_estimates}, while numerics suggest that $\theta$ is roughly $3/4$ in three dimensions. Note that our simulations yield a resistance exponent fairly close to the cut point exponent in two dimensions; it would be interesting to investigate whether these two exponents are equal or not.
 
\medskip
 
\noindent \textbf{Long-range percolation.}
For the Poisson process $\cX$ with intensity $\beta|x-y|^{-2}$, the block graph $G_\Z$ is the one-dimensional long-range percolation model in which distinct non-adjacent integers $i$ and $j$ are connected independently with probability \[1-\exp\biggl(-\beta\int_i^{i+1}\!\!\int_j^{j+1}\frac{dy\,dx}{|x-y|^{2}}\biggr).\] The existence of the exponent $\theta=\theta(\beta)$ was established for this model in \cite{ding2013distances,baumler2023distances}, while the existence of the effective resistance exponent was established in \cite{ding2025growth}. Both results establish the existence of the exponent in a stronger sense than our results and also apply to self-similar long-range percolation in all dimensions $d\geq 1$; see also \cite{ding2023uniqueness,ding2025uniqueness} for results on the scaling limit. Ding, Fan, and Huang \cite{ding2025polynomial} proved a polynomial lower bound on the effective resistance in this model, which implies that $\theta_\mathrm{res}(\beta)>0$ for every $\beta>0$. The calculation in Section~\ref{subsec:counterexamples_to_stronger_statements} shows that the expected number of cut points is of order $n^{1-\beta}$ when $\beta<1$.

\medskip
 
\noindent \textbf{Random trees and looptrees.}
For the CRT we have $\theta=\theta_\mathrm{res}=1/2$, corresponding to the fact that the CRT has dimension $2$ \cite{CRT1,aldous1993continuum}. This can be verified directly from the Brownian coding: the cut points of $G_{[0,n]}$ are the intervals containing running minima of $B$ seen from either endpoint, consecutive cut points are joined by a single pair of $\cX$, and the degrees of $G_{[0,n]}$ have exponential tails, so that effective resistances, graph distances, and cut point counts are equivalent up to a logarithmic factor. For the stable tree with index $\alpha\in(1,2)$, which has dimension $\alpha/(\alpha-1)$ \cite{duquesne2002random,duquesne2005hausdorff}, we expect that $\theta=1-1/\alpha$. Hypothesis (L) fails in this case, so that Theorem~\ref{thm:main_discrete_other_metrics} does not apply, but we expect that $\theta_\mathrm{res}$ is well-defined and equal to $\theta$ as it would be if our block graphs were actually trees. For the $\alpha$-stable looptree, which has dimension $\alpha$ \cite{curien2014random}, we expect that $\theta=1/\alpha$, in agreement with Figure~\ref{fig:numerical_estimates}(d). Since the effective resistance metric on stable looptrees is equivalent to their geodesic metric \cite{archer2021brownian}, we also expect that $\theta_\mathrm{res}=1/\alpha$ in this case.
 
\medskip
 
\noindent \textbf{Liouville quantum gravity.}
In the mating of trees approach to Liouville quantum gravity, the specific discretization considered in Theorem~\ref{thm:main_discrete} is known as the mated-CRT map.
For this model, a stronger version of Theorem~\ref{thm:main_discrete} was established by Gwynne, Holden, and Sun \cite{gwynne2019distance}, and Ding and Gwynne \cite{ding2020fractal} showed that the exponent is equal to $1/d_\gamma$, where $d_\gamma$ is the fractal dimension of $\gamma$-LQG. This dimension is computed only in the ``pure gravity'' case $\gamma=\sqrt{8/3}$, which corresponds to the scaling limit of uniform random planar maps, where $\theta=1/4$ (corresponding to the four dimensional volume growth of such maps \cite{UIPT2}). See e.g. \cite{gwynne2019distance,ding2020fractal} for bounds on $d_\gamma$ for other values of $\gamma$. On the other hand, the effective resistance on the mated-CRT map grows only logarithmically \cite{gwynne2021random}, so that $\theta_\mathrm{res}=0$ for every $\gamma\in(0,2)$.

\section{Proof of Theorems~\ref{thm:main} and \ref{thm:main_discrete}}
\label{sec:proof_main}

In this section we prove our main theorem concerning the metric exponent, Theorem~\ref{thm:main}. 
Let $\mathcal{X} \subseteq \{\{x,y\} \subset \R : x \neq y\}$ be our random set of pairs as in the statement of the theorem.
We will proceed by studying the following random variables, which are better behaved than those treated in the main theorems.
%
%
%
For each $t>0$ and $n \in \N$ define 
	\begin{align*}
		\gamma(t) = \sup_{z \in \left[0,t\right)} d^{(1)}_{\left[0,z\right]}(0,z) \qquad \text{ and } \qquad \Gamma(n) = \max_{1\leq m \leq n} \Bigl[D_{[0,m]}(0,m-1)+1\Bigr],
	\end{align*}
    with $\Gamma(0)=0$. (We add $1$ so that $\Gamma$ counts vertices rather than edges, which makes it additive under concatenation.)
    Note that $\gamma$ and $\Gamma$ are both increasing and satisfy $\gamma(t+s)\leq \gamma(t)+s$ and $\Gamma(n+1)\leq \Gamma(n)+1$, implying in particular that $\gamma$ is continuous in $t$. Note that the diameter of $[0,t]$ under the metric $d^{(1)}_{[0,t]}$ is bounded by $2\gamma(t)$ and, similarly, that the diameter of $G_{[0,n]}$ is bounded by $2\Gamma(n)$.
Our main theorem concerning these quantities does not require the full scale-invariance assumption (S), but works under the following weaker one-sided assumption:
\begin{itemize}\item[(D)] \textbf{Domination under scaling:} 
The random set $\delta \cX$ stochastically dominates $\cX$ (with respect to inclusion) for every $0<\delta\leq 1$.
\end{itemize}
(For example, if $\cX$ is a Poisson process on $\Omega$ with intensity $\beta |x-y|^{-s}$ then $\cX$ satisfies (S) when $s=2$ and satisfies (D) when $s\leq 2$.)
 Unlike our main theorem, the following theorem also establishes tightness of $\gamma$ and $\Gamma$ after rescaling, meaning that these quantities have a single order that they take with high probability (which is not true for the diameter as shown in Section~\ref{subsec:counterexamples_to_stronger_statements}).

\begin{theorem}
\label{thm:main_gamma}
Let $\mathcal{X}$ be a random closed subset of $\left\{ \{x,y\}\subset \R : x \neq y \right\}$ satisfying \emph{(I)}, \emph{(T)}, and \emph{(D)}. There exists $\theta\in [0,1]$ and a function $\sigma: (0,\infty) \to (0,\infty)$ such that $\sigma(t)=t^{\theta +o(1)}$ as $t\to\infty$ such that the families of random variables
\[
\Biggl\{
\frac{1}{\sigma(t)} \gamma(t) : t\geq 1
\Biggr\}
\qquad \text{and} \qquad
\Bigl\{ \frac{1}{\sigma(n)} \Gamma(n) : n \geq 1\Bigr\}
\]
are tight in $(0,\infty)$. Moreover, if $\mathcal{X} \neq \varnothing$ with positive probability then $\theta<1$.
\end{theorem}

\subsection{Equivalence of different formulations}

We begin by proving that the metrics from Theorem~\ref{thm:main} and Theorem~\ref{thm:main_discrete} are equivalent up-to-constants in an appropriate sense.

\begin{lemma}
\label{lem:comparison}
For each $\cX \subseteq \Omega$ we have that
\[
D_{[0,m]}(0,m-1) \leq 2 d_{[0,m]}^{(1)}(0,z)+1 \qquad \text{ and } \qquad d_{[0,m]}^{(1)}(0,z) \leq 2 D_{[0,m]}(0,m-1)+1
\]
 for each integer $m\geq 1$, and each $z\in [m-1,m]$.
\end{lemma}

	\begin{proof}[Proof of Lemma~\ref{lem:comparison}]
	Write $d=d^{(1)}_{[0,m]}$ and $D=D_{[0,m]}$, and for $x\in [0,m]$ let $c(x)=\min\{\lfloor x\rfloor, m-1\}$, so that $x\in [c(x),c(x)+1]$. Given points $u_0=0,v_0,\ldots,u_n,v_n=z$ as in the definition of $d(0,z)$, the vertices $c(v_i)$ and $c(u_{i+1})$ are equal or adjacent in $G_{[0,m]}$, while $c(u_i)$ and $c(v_i)$ are joined by a path of nearest-neighbour edges with at most $|u_i-v_i|+2$ vertices. Concatenating these paths gives a path from $0$ to $m-1$ with at most $\sum_{i=0}^n(|u_i-v_i|+2)\leq 2(n+\sum_{i=0}^n|u_i-v_i|)+2$ vertices, and taking the infimum yields $D(0,m-1)\leq 2d(0,z)+1$.
	Conversely, let $(0=x_0,x_1,\ldots,x_{\ell}=m-1)$ be a geodesic in $G_{[0,m]}$ for some $m\geq 2$, the case $m=1$ being trivial. For each $1\leq i\leq \ell$ there exist $(u_i,v_i) \in [x_{i-1},x_{i-1}+1] \times [x_{i},x_{i}+1]$ with $d(u_i,v_i)\leq 1$ (taking $u_i=v_i$ if $|x_i-x_{i-1}|=1$). Since $v_i$ and $u_{i+1}$ lie in a common unit cell, as do $0$ and $u_1$ and $v_\ell$ and $z$, we have that
		\[
			d(0,z) \leq
			d(0,u_1) + \sum_{i=1}^{\ell} d(u_i, v_i) + \sum_{i=1}^{\ell-1} d(v_i, u_{i+1}) + d(v_\ell, z)
			\leq 2\ell+1 = 2D(0,m-1)+1
		\]
		as claimed.
	\end{proof}

\subsection{The NBU property}

In this section we prove a tail estimate for the random variables $\Gamma(n)$ and $\gamma(t)$ that will play a crucial role in the proof of our main quasisubmultiplicative estimate in the next section.
    Recall that a non-negative random variable $X$ is said to have the \textbf{new better than used} (NBU) property if $\p(X \geq a+b) \leq \p(X \geq a) \p(X \geq b)$ for every $a,b\geq 0$. We say that $X$ has the integer NBU property if it is integer-valued and satisfies $\p(X \geq a+b) \leq \p(X \geq a) \p(X \geq b)$ for all integers $a,b\geq 0$.
	
	\begin{lemma}\label{lem:NBU}
		If $\cX$ satisfies \emph{(I)} and \emph{(T)} then $\gamma(t)$ has the NBU property for every $t\geq 0$ and $\Gamma(n)$ has the integer NBU property for every integer $n\geq 0$.
	\end{lemma}

In order to prove this lemma we first prove a strong Markov property for $\cX$.
For each $s\geq 0$ let $\cF_s$ be the $\sigma$-algebra generated by the pairs in $\cX$ with both endpoints in $[0,s]$, and let $\cX_s=\{\{x-s,y-s\} : \{x,y\} \in \cX,\, x,y\geq s\}$. We say that a random variable $\tau$ taking values in $[0,\infty]$ is a \textbf{stopping time} for the filtration $(\cF_s)_{s\geq 0}$ if $\{\tau \leq s\}$ is measurable with respect to $\cF_s$ for every $s\geq 0$.

\begin{lemma}[Strong Markov property]
    \label{lem:SMP}
If $\cX$ satisfies \emph{(I)} and \emph{(T)} and $\tau$ is a 
stopping time taking values in a countable set $S$, then on the event $\{\tau<\infty\}$ the random set $\cX_\tau$ is independent of $\tau$ and has the same law as $\cX_0$.
\end{lemma}

\begin{proof}[Proof of Lemma~\ref{lem:SMP}]
Sum over the possible values $s\in S$ of $\tau$, using that $\{\tau=s\}\in\cF_s$ (since $S$ is countable) and that $\cX_s$ is independent of $\cF_s$ by (I) and has the law of $\cX_0$ by (T) for each $s\in S$.
\end{proof}

\begin{remark}
\emph{
Lemma~\ref{lem:SMP} is false without the assumption that $\tau$ is countable-valued. Indeed, take $\cP$ to be a Poisson process on $\R$ and take $\cX$ to be the set of pairs $\{\{p-1,p\}:p\in \cP\}\cup \{\{p,p+1\}:p\in \cP\}$. It is easily verified that $\cX$ satisfies (I) and (T), but  if we define $\tau$ to be the infimal value of $t$ such that the restriction of $\cX$ to $[0,t]$ is non-trivial then $\cX_\tau$ does not have the same distribution as $\cX_0$ as it includes the pair $\{0,1\}$ with positive probability.
}
\end{remark}

\begin{proof}[Proof of Lemma~\ref{lem:NBU}]
We prove the claim for $\Gamma(n)$, the proof for $\gamma(t)$ being similar (one simply applies the same argument to the stopping times $2^{-j}\lceil 2^j T_a\rceil$, where $T_a=\inf\{s:\gamma(s)\geq a\}$, then takes $j\to\infty$). For each integer $a \geq 0$ let $N_a=\inf\{m\geq 0 : \Gamma(m)\geq a\}$. It suffices to prove that $N_{a+b}$ stochastically dominates the sum of two independent copies of $N_a$ and $N_b$. Indeed, the statement of the lemma is equivalent to $N_{a+b}$ stochastically dominating the \emph{maximum} of two independent copies of $N_a$ and $N_b$, which is weaker than this claim. To prove this claim, observe that if $N_b'$
is the smallest value of $n$ such that $D_{[N_a,N_a+n]}(N_a,N_a+n-1)+1 \geq b$ then $N_b'$ is independent of $N_a$ and has the same law as $N_b$ by Lemma~\ref{lem:SMP}, and the claim follows by noting that $N_{a+b} \geq N_a + N_b'$ by the triangle inequality (counting vertices) and the fact that $D_{[a,b]}(x,y)$ is decreasing in the choice of interval $[a,b]$.
\end{proof}

It follows from Lemma~\ref{lem:NBU} that if we define $M(n)=\inf\{\ell\in \N : \p(\Gamma(n) \geq \ell) \leq 1/2\}$ to be the median of $\Gamma(n)$ then
\begin{equation}
\label{eq:NBU_tail_Gamma}
\p\Bigl( \Gamma(n) \geq k M(n) \Bigr) \leq 2^{-k}
\qquad\text{and}\qquad
 \p\Bigl(\Gamma(n) \geq \Bigl\lfloor \frac{M(n)-1}{k}\Bigr\rfloor\Bigr) \geq 2^{-1/k}
\end{equation}
for every integer $k\geq 1$, with similar tail estimates holding for $\gamma(t)$.

\subsection{Quasisubmultiplicativity}

We now prove the following quasisubmultiplicative property of $\E\Gamma(n)$ which will suffice to prove the existence of the exponent. A more careful version of the analysis will be required to prove sublinearity of the exponent for non-trivial $\cX$.

	\begin{proposition}\label{prop:improved submult}
	If $\cX$ satisfies \emph{(I)}, \emph{(T)}, and \emph{(D)} then 
		\begin{equation*}
			\E \Gamma(mn) \leq  \frac{36e}{\log 2}  \left(1 + \log \frac{n}{\E \Gamma(n)} \right) \E \Gamma(m) \E \Gamma (n) 
		\end{equation*}
		for every $m,n \geq 2$.
	\end{proposition}
	
In order to prove this proposition, we first prove the following general lemma about random variables with exponential tails.

    \begin{lemma}
    \label{lem:order_statistics}
Let $c>0$ and $C\geq 1$ be positive constants and suppose that $X_1,\ldots,X_m$ are (not necessarily independent or identically distributed) non-negative random variables each satisfying the tail bound $\p(X_i \geq x) \leq Ce^{-cx}$ for all $x\geq 0$. If $A \subseteq \{1,\ldots,m\}$ is a random set of indices defined on the same probability space as $(X_1,\ldots,X_m)$ then
\[
\E \sum_{i\in A} X_i \leq \frac{e}{c}\E|A| \left(1+ \log \frac{Cm}{\E|A|}\right).
\]
In particular, if we rearrange the sequence $X_1,\ldots,X_m$ in decreasing order $X_{(1)}\geq X_{(2)} \geq \cdots \geq X_{(m)}$ then
\[
\E \sum_{i=1}^N X_{(i)} \leq \frac{e}{c} \E N\left(1+\log \frac{Cm}{\E N}\right)
\]
for every random variable $N$ taking values in $\{1,\ldots,m\}$ and defined on the same probability space as $(X_1,\ldots,X_m)$.
    \end{lemma}

    \begin{proof}[Proof of Lemma~\ref{lem:order_statistics}]
Let $\pi:\{1,\ldots,m\}\to\{1,\ldots,m\}$ be a uniform random permutation independent of $(X_1,\ldots,X_m)$ and $A$, and observe that 
$\p(X_{\pi(1)}\geq x) \leq Ce^{-cx}$ and that
\begin{align*}
\E \sum_{i\in A} X_i = \E \sum_{i:\pi(i)\in A} X_{\pi(i)} = m \E X_{\pi(1)} \mathbbm{1}(\pi(1)\in A).
\end{align*}
We can also compute that
\[
\E X_{\pi(1)}^p = \int_0^\infty px^{p-1}\p(X_{\pi(1)} \geq x) dx \leq C  \int_0^\infty px^{p-1} e^{-cx} dx = C c^{-p} \Gamma(p+1)
\]
for every $p>0$, 
where $\Gamma$ denotes the $\Gamma$ function, 
and hence by 
H\"older's inequality that
\begin{align*}
\E \sum_{i\in A} X_i  \leq m \p(\pi(1)\in A)^{(p-1)/p}
(\E X_{\pi(1)}^p)^{1/p} \leq m C^{1/p} c^{-1} (\E|A|/m)^{(p-1)/p} \Gamma(p+1)^{1/p}
\end{align*}
for every $p\geq 1$. Taking $p=1\vee \log Cm/\E|A|$ and using that $\Gamma(p+1)^{1/p}\leq p$ for every $p\geq 1$ yields the claim.
    \end{proof}

	\begin{proof}[Proof of Proposition~\ref{prop:improved submult}]
		For each $m,n\geq 1$, let $G_{[0,n]}^m$ be the graph defined like $G_{[0,n]}$ but using the re-scaled process $\frac{1}{m}\cX$ instead of $\cX$. In other words, this graph has vertex set $\left\{0,\ldots,n-1\right\}$ and has two points $i$ and $j$ of this interval adjacent if either $|i-j|=1$ or $|i-j|>1$ and there exists $\{x,y\}\in \cX$ such that $x \in [mi,m(i+1)]$ and $y\in [mj,m(j+1)]$. Writing $D^m_{[0,n]}$ for the graph metric on $G_{[0,n]}^m$, the assumption (D) implies that $\Gamma^m(n)= \max_{1\leq \ell \leq n} [D^m_{[0,\ell]}(0,\ell-1)+1]$ is stochastically dominated by $\Gamma(n)$, while we also have that $\Gamma^m(n)$ satisfies the integer NBU property since $\frac{1}{m}\cX$ satisfies (I) and (T).

        Now, define the random variables $X_0,\ldots,X_{n-1}$ by
		\begin{equation*}
			X_i = \max_{1\leq \ell \leq m} \Bigl[D_{\left[im,im+\ell \right]}(im,im+\ell-1)+1\Bigr],
		\end{equation*} 
		which are independent and all distributed as $\Gamma(m)$ by (I) and (T), and observe that the diameter of the graph $G_{[im,(i+1)m]}$ is at most $2 X_i$ for each $0\leq i \leq n-1$.
We claim that if $X_{(1)}\geq X_{(2)} \geq \cdots \geq X_{(n)}$ is the decreasing rearrangement of the $X_i$ then 
\begin{equation}
\label{eq:submult_pointwise}
	\Gamma(mn) \leq 2 \sum_{i=1}^{\min\{\Gamma^m(n)+2,n\}} X_{(i)}.
\end{equation}
To prove this, we fix $p=\ell m + r$ with $0\leq \ell \leq n-1$ and $1\leq r \leq m$. If $\ell=0$ then $D_{[0,p]}(0,p-1)+1\leq X_0$, so suppose that $\ell\geq 1$. We can construct a path from $0$ to $p-1$ in $G_{[0,p]}$ as follows: Start with a geodesic $x_0,\ldots,x_a$ from $0$ to $\ell-1$ in $G^m_{[0,\ell]}$. For each increment $(x_j,x_{j+1})$ in this path, there are points $z_j^+ \in [mx_j,m(x_{j}+1)]$ and $z_{j+1}^- \in [mx_{j+1},m(x_{j+1}+1)]$ such that either $z_j^+=z_{j+1}^-$ or $\{z_j^+,z_{j+1}^-\}\in \cX$ (the first case is to allow for $|x_j-x_{j+1}|=1$). Setting $z_0^-=0$ and $z_a^+=\ell m$, for each $0\leq j\leq a $ there is a path with at most $2X_{x_j}$ vertices in $G_{[mx_j,m(x_j+1)]}$ between cells containing $z_j^-$ and $z_j^+$, and consecutive such paths can be joined by at most one edge. We can concatenate all of these paths with the edge $\{\ell m-1,\ell m\}$ and a geodesic from $\ell m$ to $p-1$ in $G_{[\ell m,p]}$, which has at most $X_\ell$ vertices, to obtain a path from $0$ to $p-1$ in $G_{[0,p]}$. Since the geodesic $x_0,\ldots,x_a$ visits each point in $\{0,\ldots,\ell-1\}$ at most once and $a+1\leq \Gamma^m(n)$, this path has at most $2 \sum_{i=1}^{a+2} X_{(i)}\leq 2 \sum_{i=1}^{\min\{\Gamma^m(n)+2,n\}} X_{(i)}$ vertices as claimed.

		Let $M(m)=\min\{k \geq 1: \p(\Gamma(m) \geq k) \leq 1/2\}$ be the median of $\Gamma(m)$, which satisfies $\E \Gamma(m) \geq \frac{1}{2} (M(m)-1)$ by definition. The estimate \eqref{eq:NBU_tail_Gamma} implies that each of the random variables $X_i$ satisfies 
		\[
			\p(X_i \geq k M(m)) \leq 2^{-k}
		\]
		for every $k\geq 1$, so that they also satisfy 
		\[
			\p(X_i \geq x) \leq \exp\left[1- \log 2 \frac{x}{M(m)}\right]
		\]
		for every real $x>0$.
		As such, it follows from \eqref{eq:submult_pointwise} and Lemma~\ref{lem:order_statistics} (applied with $C=e$, $c=\log 2/M(m)$, and $N=\min\{\Gamma^m(n)+2,n\}$) that
		\begin{align*}
			\E \Gamma(mn) 
			&\leq 2e \frac{\E (\Gamma^m(n)+2) M(m)}{\log 2} \left(2 + \log \frac{n}{\E (\Gamma^m(n)+2)} \right)
			\\
			&
			\leq 36 e \frac{\E \Gamma^m(n) \E \Gamma(m)}{\log 2} \left(1 + \log \frac{n}{\E \Gamma^m(n)} \right) \leq 36 e \frac{\E \Gamma(n) \E \Gamma(m)}{\log 2} \left(1 + \log \frac{n}{\E \Gamma(n)} \right)
		\end{align*}
		as claimed, where we used that $1\leq \Gamma^m(n)\leq n$ and $M(m) \leq 2 \E\Gamma(m)+1 \leq 3 \E\Gamma(m)$ to simplify the expression and used that $x \mapsto x(c+\log (n/x))$ is increasing on $(0,e^{c-1}n]$ for $c\in\{1,2\}$ to replace $\E N$ by $\E(\Gamma^m(n)+2)$ and to apply the stochastic domination estimate $\E \Gamma^m(n) \leq \E \Gamma(n)$ in the final inequality.\qedhere
	\end{proof}	
	
We now apply Proposition~\ref{prop:improved submult} to prove the existence of the exponent from Theorem~\ref{thm:main_gamma}.

	\begin{corollary}
	\label{cor:first_exponent}
		If $\cX$ satisfies \emph{(I)}, \emph{(T)}, and \emph{(D)} then there exists $\theta \in \left[0,1\right]$ such that
			$\E \Gamma(n) = n^{\theta + o(1)}$ as $n\to \infty$.
	\end{corollary}
	
	\begin{proof}[Proof of Corollary~\ref{cor:first_exponent}]
	Write $f(n)=\E \Gamma(n)$.
		Proposition~\ref{prop:improved submult} directly implies that the weaker form of the submultiplicativity
		\begin{align*}
			f(mn) \leq C \log(m) f(m) f(n)
		\end{align*}
		holds for all $m,n \geq 2$ and some constant $C< \infty$. This implies that the sequence $a_k=\log_2 f(2^k)$ satisfies the subadditive inequality
		\[
			a_{m+n} \leq a_m + a_n + \log_2 C m.
		\]
		Since the error term grows only logarithmically, the de Bruijn--Erd\accent"7D os generalization of Fekete's lemma \cite{de1952some} implies that $\theta=\lim_{k\to \infty} a_k/k$ is well-defined, so that $f(2^k)=2^{k(\theta+o(1))}$ as $k\to \infty$. This implies that $f(n)=n^{\theta+o(1)}$ as $n\to\infty$ since $f(n)$ is monotone in $n$, and the fact that $\theta \in [0,1]$ follows from the trivial bounds $1\leq f(n) \leq n$.
	\end{proof}

\subsection{Sublinearity of the diameter for non-trivial quotients}
	
In this section we prove that if $\cX$ is non-trivial then $\theta<1$. This will follow from the quasisubmultiplicativity property of Proposition~\ref{prop:improved submult} once we prove the following sublinearity estimate. 

	\begin{lemma}\label{lem:o(n)}
		Suppose that $\cX$ satisfies \emph{(I)}, \emph{(T)}, and \emph{(D)}. If $\p(\cX \neq \emptyset ) > 0$,
        then $\E \Gamma(n) = o(n)$.
	\end{lemma}

	The proof of this lemma is similar to the proof of \cite[Proposition 2]{curien2020geometric}.
	
	\begin{proof}[Proof of Lemma~\ref{lem:o(n)}]
		Kingman's subadditive ergodic theorem, applied to the subadditive family $(D_{[m,n]}(m,n-1)+1)_{0\leq m<n}$, implies that the random variable
		\begin{equation}\label{X_n limit}
			Z=\lim_{n \to \infty} \frac{1}{n} D_{[0,n]}(0,n-1) 
		\end{equation}
		is well-defined, and it follows by bounded convergence that $\frac{1}{n}\E \Gamma(n) \to \E Z$. As such, it suffices to prove that $Z=0$ almost surely.
		Let $\eps \in (0,0.5)$ be arbitrary. Since $\cX$ contains a pair in some set of the form $[a,a+h]\times[b,b+h]$ with $0<h\leq \eps(b-a)$ with positive probability, the stochastic monotonicity of $\frac{1}{n} \cX$ (property (D)) and translation invariance (property (T)) imply that there exists $n_0=n_0(\eps)$ and $\delta=\delta(\eps)>0$ such that 
		\begin{equation*}
		\p(A_n):=	\p \left(\exists x \in \left[\eps n, 2\eps n\right], y \in \left[n-\eps n, n\right] : \{x,y\} \in \cX\right) \geq \delta
		\end{equation*}
        for all $n\geq n_0$. 
Since the events $A_{n_1}$ and $A_{n_2}$ are independent when $\eps n_2 \geq n_1$, infinitely many of these events occur almost surely. Since $D_{[0,n]}(0,n-1) \leq 3\eps n+3$ whenever $A_n$ holds, this implies that $Z\leq 3\eps$ almost surely, and the claim follows since $\eps>0$ was arbitrary. \qedhere
	\end{proof}

	\begin{proof}[Proof of Theorem~\ref{thm:main_gamma}]
	Lemma~\ref{lem:comparison} implies that $\gamma(t)\leq 2\Gamma(\lceil t\rceil)$ and $\Gamma(\lceil t \rceil) \leq 2\gamma(t)+3$ for every $t>0$. (For the first inequality, note that if $z\in[m-1,m)$ for some integer $m\geq 2$ then $d^{(1)}_{[0,z]}(0,z)\leq d^{(1)}_{[0,m-1]}(0,m-1)+1\leq 2\Gamma(m-1)$ by Lemma~\ref{lem:comparison} applied with $m-1$ in place of $m$, while $d^{(1)}_{[0,z]}(0,z)\leq z< 1$ if $z<1$.) Given this, Lemma~\ref{lem:NBU}, and Corollary~\ref{cor:first_exponent}, it remains only to prove that $\theta<1$ when $\p(\cX \neq \varnothing)>0$.
		Let $C=72e/\log 2$ be the constant from Proposition~\ref{prop:improved submult}, and let $\eps \in (0,0.1)$ small enough that $C(1+\log(\frac{1}{\eps})) \eps \leq 1/2$. It follows from Lemma~\ref{lem:o(n)} that there exists $n_0\geq 2$ such that $\E \Gamma(n_0) \leq \eps n_0$, and it follows inductively from Proposition~\ref{prop:improved submult} (using that $x(1+\log(1/x))$ is increasing on $(0,1]$) that $\E \Gamma(n_0^k) \leq \eps 2^{-k+1} n_0^k$ for every $k\geq 1$. This implies that $\theta <1$ as claimed. \qedhere
	\end{proof}

\subsection{The diameter}
\label{subsec:the_diameter}


To prove Theorems~\ref{thm:main} and \ref{thm:main_discrete}, it remains to prove that $\E d_{[0,t]}^{(1)}(0,t) \geq t^{\theta-o(1)}$. We write $d_{[0,t]}=d^{(1)}_{[0,t]}$ to lighten notation, and begin by noting that $\E d_{[0,t]}(0,t)$ is increasing in $t$.

\begin{lemma}
\label{lem:domination_time}
If $\cX$ satisfies \emph{(S)} then $d_{[0,t]}(0,t)$ stochastically dominates $d_{[0,s]}(0,s)$ for every $0\leq s \leq t$.
\end{lemma}

\begin{proof}[Proof of Lemma~\ref{lem:domination_time}]
The distance $d_{[0,s]}(0,s)=d_{[0,s]}^{(1)}(0,s)$ has the same distribution as $\frac{s}{t} d^{(t/s)}_{[0,t]}(0,t)$, which is bounded by $d^{(1)}_{[0,t]}(0,t)$ when $t\geq s$.
\end{proof}

We next show that $\E d_{[0,t]}(0,t)$ and $\E \Gamma(t)$ can be related for certain well-chosen choices of $t$. Given $\eta>0$, we say that $t$ is \textbf{$\eta$-regular} if 
\[
	t^{-\eta} \E \gamma(t) \geq \eta \sup_{s\leq t} s^{-\eta} \E \gamma(s).
\]
The following (trivial) lemma uses only that $\E \gamma(t)$ has a well-defined rate of polynomial growth.

\begin{lemma}
\label{lem:regular_times}
Let $\cX$ satisfy \emph{(I)}, \emph{(T)}, and \emph{(D)} and let $\theta$ be as in Theorem~\ref{thm:main_gamma}. If $\theta>0$ then for each $\eps>0$ and $0<\eta<\theta$ there exists $t_0=t_0(\eps,\eta)<\infty$ such that there exists an $\eta$-regular $s$ in $[t^{1-\eps},t]$ for every $t\geq t_0$.
\end{lemma}

\begin{proof}[Proof of Lemma~\ref{lem:regular_times}]
Note that for $s$ to be $\eta$-regular it suffices for $s$ to be a record time for $s^{-\eta} \E \gamma(s)$ in the sense that $s^{-\eta} \E \gamma(s) = \sup_{s'\leq s} (s')^{-\eta} \E \gamma(s')$.  
By (the proof of) Theorem~\ref{thm:main_gamma}, for each $\delta>0$ there exists $t_1$ such that $t^{\theta-\delta}\leq \E \gamma(t) \leq t^{\theta+\delta}$ for every $t\geq t_1$. Thus, if $\theta-\delta-\eta>0$ and $t$ is large then $s^{-\eta} \E \gamma(s) \leq \max\{t_1,s^{\theta+\delta-\eta}\} \leq t^{(1-\eps)(\theta+\delta-\eta)}$ for every $s\leq t^{1-\eps}$ and $t^{-\eta} \E \gamma(t) \geq t^{\theta-\delta-\eta}$, so that $s^{-\eta}\E\gamma(s)$ must have a record between $t^{1-\eps}$ and $t$ when $t$ is large and $\theta-\delta-\eta > (1-\eps)(\theta+\delta-\eta)$. The claim follows since $\eta<\theta$ and $\eps,\delta>0$ were arbitrary.
\end{proof}

We next show that $d_{[0,t]}(0,t)$ has good probability to be of order $\E \gamma(t)$ when $t$ is $\eta$-regular.

\begin{lemma}
\label{lem:gamma_to_diam}
Let $\cX$ satisfy \emph{(I)}, \emph{(T)}, and \emph{(S)}. For each $\eta>0$ there exists a positive constant $c(\eta)$ such that if $t\geq 1$ is $\eta$-regular then
\[
\p\Bigl(d_{[0,t]}(0,t) \geq \frac{1}{4}\E \gamma(t)\Bigr) \geq c(\eta). 
\]
\end{lemma}

\begin{proof}[Proof of Lemma~\ref{lem:gamma_to_diam}]
Let $M(t)=\inf\{r\geq 0: \p(\gamma(t) \geq r)\leq 1/2\}$ be the median of $\gamma(t)$, noting that
\begin{multline*}
	\frac{1}{2}M(t) \leq \E \gamma(t) = \int_0^\infty \p(\gamma(t) > x) \md x \\\leq M(t) \sum_{k \geq 0} \p(\gamma(t)> k M(t)) \leq M(t) \sum_{k=0}^\infty 2^{-k}= 2M(t)
\end{multline*}
by the NBU property. Write $f(z)=d_{[0,z]}(0,z)$ for each $z\geq 0$, and note that $f(z')\leq f(z)+z'-z$ for every $0\leq z\leq z'$.
Suppose that 
$\gamma(t) \geq M(t)$, let $t_0\in [0,t)$ be minimal such that $f(t_0)\geq \frac{3}{4}\gamma(t)$, and let $0\leq s_0 < t_0$ be minimal such that 
\[
	\frac{|\{s \in [s_0,t_0] : f(s) \geq \frac{1}{2} \gamma(t)\}|}{|s_0-t_0|} \geq \frac{1}{2},
\]
where the quantity in the numerator is the Lebesgue measure of the relevant set.
Note that $s_0<t_0$ since $f(s)\geq f(t_0)-(t_0-s)$, and that if $s_0\neq 0$ then the set $B=\{s \in [s_0,t_0] : f(s) < \frac{1}{2} \gamma(t)\}$ has measure $\frac{1}{2}|s_0-t_0|$, while $d_{[s,t_0]}(s,t_0)\geq f(t_0)-f(s) \geq \frac{1}{4}\gamma(t)$ for every $s\in B$ by the triangle inequality.
For each 
$k\geq 1$ let $\cE(k)$ be the event that 
 $2^{-k} t \leq |s_0-t_0| < 2^{-k+1}t$. On the event $\cE(k) \cap \{\gamma(t) \geq  M(t)\} \cap \{s_0 \neq 0\}$, if we pick a point $X$ uniformly at random from the interval $[0,t]$ then with probability at least $2^{-k-1}$ it belongs to the set $B$ and therefore satisfies
\[\sup_{0\leq r < 2^{-k+2}t} d_{[X,X+r]}(X,X+r) \geq d_{[X,t_0]}(X,t_0) \geq \frac{1}{4} \gamma(t) \geq \frac{1}{4} M(t),\]
where the left hand side has the same law as $\gamma(2^{-k+2}t)$ by (T).
Meanwhile if $\cE(k) \cap \{\gamma(t) \geq M(t)\} \cap \{s_0 = 0\}$ holds then we must have that $t_0< 2^{-k +1} t$ and hence that $\gamma(2^{-k+1} t) \geq f(t_0)\geq \frac{3}{4} M(t)$. Taking a union bound over these two possibilites yields that
\begin{multline*}
	\p(\cE(k)\cap \{\gamma(t)\geq  M(t)\}) \leq 2^{k+1} \p\Bigl(\gamma(2^{-k+2}t) \geq \frac{1}{4}M(t)\Bigr)+\p\Bigl(\gamma(2^{-k+1}t) \geq  \frac{3}{4}M(t)\Bigr) 
	\\\leq 2^{k+2} \p\Bigl(\gamma(2^{-k+2}t) \geq \frac{1}{4}M(t)\Bigr).
\end{multline*}
Since $t$ is $\eta$-regular, the NBU property bounds this probability for $k\geq 2$ by
\begin{multline*}
	-\log_2 \p(\cE(k)\cap \{\gamma(t)\geq  M(t)\}) \geq -k-3 + \frac{ M(t)}{4 M(2^{-k+2} t)} 
	\\
	\geq -k-3 + \frac{ \E \gamma(t)}{16 \E \gamma(2^{-k+2} t)} 
	\geq -k-3 + \frac{\eta}{16} 2^{\eta (k-2)}.
\end{multline*}
It follows that there exists $k_0=k_0(\eta)\geq 2$ such that
\[
	\sum_{k=k_0+1}^\infty \p(\cE(k)\cap \{\gamma(t)\geq  M(t)\}) \leq \frac{1}{4}
\]
and hence that
\[
	\p(\gamma(t) \geq M(t) \text{ and } |s_0-t_0|\geq 2^{-k_0} t) \geq \frac{1}{4}.
\]
On this event, if we pick $X$ uniformly at random from $[0,t]$ then it satisfies $d_{[0,X]}(0,X) \geq \frac{1}{2}\gamma(t)\geq \frac{1}{2}M(t)$  with probability at least $2^{-k_0-1}$. As such, Fubini's theorem implies that there exists $s\in [0,t]$ such that $\p(d_{[0,s]}(0,s)\geq \frac{1}{2}M(t))\geq 2^{-k_0-3}$, and the claim follows from Lemma~\ref{lem:domination_time} and the inequality $M(t) \geq \E \gamma(t)/2$.
\end{proof}

\begin{proof}[Proof of Theorems~\ref{thm:main} and \ref{thm:main_discrete}]
Let $\theta\in [0,1]$ be as in Theorem~\ref{thm:main_gamma}.
We have trivially that 
\[\E d^{(1)}_{[0,t]}(0,t) \leq \E \sup_{x,y\in [0,t]} d^{(1)}_{[0,t]}(x,y) \leq 2\E \gamma(t) = t^{\theta+o(1)}.\] If $\theta=0$ then the complementary lower bound holds trivially. Otherwise, if we fix $0<\eta<\theta$ and let $s(t)=\max\{s\leq t: s$ is $\eta$-regular$\}$ for each $t>0$ then $s(t)=t^{1-o(1)}$ by Lemma~\ref{lem:regular_times}
and 
\[\E d^{(1)}_{[0,t]}(0,t) \geq \E d^{(1)}_{[0,s(t)]}(0,s(t)) \geq \frac{c(\eta)}{4} \E \gamma(s(t)) = t^{\theta-o(1)}\] by Lemma~\ref{lem:domination_time} and Lemma~\ref{lem:gamma_to_diam}. The remaining claims follow by the scaling identity \eqref{eq:basic_scaling_identity} and Lemma~\ref{lem:comparison}, which implies that $D_{[0,n]}(0,n-1)\geq \frac{1}{2}(d^{(1)}_{[0,n]}(0,n)-1)$ for every $n\geq 1$.
\end{proof}

\section{Effective resistances}
\label{sec:effective_resistances}

In this section we prove our theorem concerning the existence of the growth exponent for the effective resistance, Theorem~\ref{thm:main_discrete_other_metrics}. We begin by introducing the \emph{vertex extremal length} (VEL) and explaining its relation to the effective resistance in Section~\ref{subsec:vertex_extremal_length}. We then adapt the argument of Section~\ref{sec:proof_main} to the VEL in Sections~\ref{subsec:generalizing_the_quasisubmultiplicative_argument} and \ref{subsec:point_to_point}; this implies that the effective resistance exponent is well-defined when (L) holds since the effective resistance and the VEL are equivalent up to a factor of the maximum degree. Finally, in Section~\ref{subsec:L_brownian} we prove that (L) holds for the trace of Brownian motion in two and three dimensions.

\subsection{Vertex extremal length}
\label{subsec:vertex_extremal_length}

We recall that the effective resistance between two points on a finite graph can be expressed as
\begin{equation}
	\mathscr{R}_\mathrm{eff}(u\leftrightarrow v) = \mathscr{R}_\mathrm{eff}(u\leftrightarrow v;G) = \inf_\p \sum_{e\in E} \p(e \in \mathscr{P})^2 
	\label{eq:resistance_def}
\end{equation}
where the infimum is taken over probability measures on paths between $u$ and $v$. Using the minimax theorem, this quantity can be expressed alternatively as
\begin{align}
\inf_\p \sum_{e\in E} \p(e \in \mathscr{P})^2 &= \left(\inf_\p \sup_{\mathbb{Q}} \E_{\p \otimes \mathbb{Q}} \left[\sum_{e\in E} \omega(e)\mathbbm{1}(e\in \mathscr{P})\right]\right)^2
\nonumber
\\&= 
\left(\sup_{\mathbb{Q}} \inf_\p  \E_{\p \otimes \mathbb{Q}} \left[\sum_{e\in E} \omega(e)\mathbbm{1}(e\in \mathscr{P})\right]\right)^2 = \sup_{\omega} d_\omega(u,v)^2,
\label{eq:resistance_duality}
\end{align}
where $\p$ is chosen from the set of probability measures on paths between $u$ and $v$ as before, $\omega$ is chosen from the set of functions $\omega:E\to [0,1]$ with $\sum_{e\in E} \omega(e)^2\leq 1$, and $\mathbb{Q}$ is chosen from the set of probability measures on such functions. Note here that given the probability measure $\p$, the optimal choice of $\mathbb{Q}$ is equal to the Dirac mass on the function
\[
	\omega(e)=\frac{\p(e\in \mathscr{P})}{\sqrt{\sum_{e\in E} \p(e\in \mathscr{P})^2}}
\]
while, given $\mathbb{Q}$, an optimal choice of $\p$ is given by the Dirac mass on some shortest path between $u$ and $v$ in the metric with edge weights $\mathbb{Q}[\omega(e)]$; taking $\mathbb{Q}$ to be a Dirac mass is always optimal since $\sum_{e\in E} \mathbb{Q}[\omega(e)]^2 \leq \sum_{e\in E} \mathbb{Q}[\omega(e)^2]$ with equality if and only if $\mathbb{Q}$ is supported on a single function $\omega$. This expression for the effective resistance is also known as the \textbf{edge extremal length} (EEL).

We will also work with the \textbf{vertex extremal length}, which was introduced by Cannon \cite{cannon1994combinatorial} and plays an important role in the theory of circle packing \cite{HeSc}. Given a finite graph $G=(V,E)$ and a pair of distinct vertices $u,v\in V$, we define the vertex extremal length (VEL) between $u$ and $v$ to be
\[
	\mathscr{V}(u\leftrightarrow v) = \mathscr{V}(u\leftrightarrow v; G) = \inf \sum_{w\in V} \p(w \in \mathscr{P})^2,
\]
where again the infimum is taken over probability measures on paths between $u$ and $v$. By the same minimax calculation above, the vertex extremal length can also be expressed as
\begin{align}
\inf_\p \sum_{w\in V} \p(w \in \mathscr{P})^2 &= \left(\inf_\p \sup_{\mathbb{Q}} \E_{\p \otimes \mathbb{Q}} \left[\sum_{w\in V} \omega(w)\mathbbm{1}(w\in \mathscr{P})\right]\right)^2
\nonumber
\\&= 
\left(\sup_{\mathbb{Q}} \inf_\p  \E_{\p \otimes \mathbb{Q}} \left[\sum_{w\in V} \omega(w)\mathbbm{1}(w\in \mathscr{P})\right]\right)^2 = \sup_{\omega} d_\omega(u,v)^2,
\label{eq:VEL_duality}
\end{align}
where $\omega$ is now a function $V\to [0,1]$ with $\sum_{w\in V}\omega(w)^2\leq 1$, $\mathbb{Q}$ is now a probability measure on functions of the same form, and $d_\omega$ is defined to be the infimum of $\sum_{w\in V} \omega(w)\mathbbm{1}(w\in \mathscr{P})$ over paths $\mathscr{P}$ between $u$ and $v$. (Note that the VEL is usually defined through this dual formulation.) Note that the random path attaining the infimum in the definition of the vertex extremal length can always be taken to be almost surely simple.


The following lemma will permit us to prove our results for the VEL instead of the effective resistance under the assumption (L).

\begin{lemma}
\label{lem:EEL_VEL}
If $G$ is a graph with maximal degree $\Delta$ then
\[
\frac{1}{2\Delta}	\mathscr{V}(u\leftrightarrow v)\leq 
	\mathscr{R}_\mathrm{eff}(u\leftrightarrow v) \leq \mathscr{V}(u\leftrightarrow v)-1
\]
for every distinct $u,v\in V$.
\end{lemma}

\begin{proof}[Proof of Lemma~\ref{lem:EEL_VEL}]
The sums $\sum_{e\in E}\p(e\in \mathscr{P})^2$ and $\sum_{w\in V}\p(w\in \mathscr{P})^2$ are equal to the expected numbers of edges and vertices shared by two independent copies of $\mathscr{P}$ respectively. The upper bound follows since the number of edges shared by two simple random paths is at most one less than the number of vertices they share.
We now prove the lower bound. If $\omega:V\to [0,1]$ has $\sum_{w\in V} \omega(w)^2\leq 1$ then the function $\omega':e\to [0,1]$ defined by $\omega'(e)= \frac{1}{\sqrt{2\Delta}}(\omega(x)+\omega(y))$ for each edge $e$ with endpoints $x$ and $y$ satisfies
\[
	\sum_{e\in E} \omega'(e)^2 \leq  \sum_{w\in V} \omega(w)^2 \leq 1
\qquad \text{and} \qquad
	d_{\omega'}(u,v) \geq \frac{1}{\sqrt{2\Delta}} d_\omega(u,v),
\]
where the second inequality holds since every vertex of a path belongs to at least one of its edges. This is easily seen to imply the claim.
\end{proof}

\begin{remark}
\emph{It is a classical fact \cite[Lemma 4.0.1]{KumagaiBook} that the effective resistance $\mathscr{R}_\mathrm{eff}(u\leftrightarrow v)$ defines a metric on any finite graph.
It can be shown that $\mathscr{V}(x\leftrightarrow y)-1$ also defines a metric on any finite graph using essentially the same argument used for the $L^p$ analogue of the edge extremal length in \cite{albin2019blocking}. We do not give the proof of this as we do not need it.}
\end{remark}

\subsection{Generalizing the quasisubmultiplicative argument}
\label{subsec:generalizing_the_quasisubmultiplicative_argument}
	
We now generalize the quasisubmultiplicative analysis of Proposition~\ref{prop:improved submult} to the VEL. For integers $a<b$ and $x,y\in \{a,\ldots,b-1\}$ we write
\[
D^\mathrm{res}_{[a,b]}(x,y) = \mathscr{R}_\mathrm{eff}(x\leftrightarrow y; G_{[a,b]}) \qquad \text{and} \qquad D^\mathrm{vel}_{[a,b]}(x,y) = \mathscr{V}(x\leftrightarrow y; G_{[a,b]}),
\]
and for each $\# \in \{\mathrm{res},\mathrm{vel}\}$ define
\[\Gamma_\#(n) = \max_{1\leq m \leq n} D^\#_{[0,m]}(0,m-1) \quad \text{for $n\geq 1$} \qquad \text{and} \qquad \Gamma_\#(0)=0.\]
The proof of Lemma~\ref{lem:NBU} goes through essentially unmodified for the VEL.

	\begin{lemma}\label{lem:NBU_generalized}
		If $\cX$ satisfies \emph{(I)} and \emph{(T)} then $\p(\Gamma_\mathrm{vel}(n)\geq a+b+1)\leq \p(\Gamma_\mathrm{vel}(n)\geq a)\p(\Gamma_\mathrm{vel}(n)\geq b)$ for every integer $n\geq 0$ and every $a,b\geq 0$.
	\end{lemma}

\begin{proof}[Proof of Lemma~\ref{lem:NBU_generalized}]
The proof is identical to that of Lemma~\ref{lem:NBU}: that proof used only that 
\[
	D^\mathrm{vel}_{[0,m]}(0,m-1) \leq D^\mathrm{vel}_{[0,k]}(0,k-1) + D^\mathrm{vel}_{[k,m]}(k,m-1)
\]
for every $0<k<m$ (which holds by concatenating independent optimal random paths in $G_{[0,k]}$ and $G_{[k,m]}$, which have disjoint vertex sets, via the edge $\{k-1,k\}$), that $D^\mathrm{vel}_{[k,m]}$ is determined by the restriction of $\cX$ to $[k,m]$, and that $\Gamma_\mathrm{vel}(m+1)\leq \Gamma_\mathrm{vel}(m)+1$. The additional $+1$ accounts for the fact that $\Gamma_\mathrm{vel}$ is not integer-valued.
\end{proof}

We next adapt the quasisubmultiplicativity estimate of Proposition~\ref{prop:improved submult}. The proof does not appear to work for the effective resistance, which is why we work with the vertex extremal length instead.

	\begin{proposition}\label{prop:improved_submult_generalized}
	There exists a universal constant $C$ such that if $\cX$ satisfies \emph{(I)}, \emph{(T)}, and \emph{(D)} then 
		\begin{equation*}
			\E \Gamma_\mathrm{vel}(mn) \leq  C  \left(1 + \log \frac{n}{\E \Gamma_\mathrm{vel}(n)} \right) \E \Gamma_\mathrm{vel}(m) \E \Gamma_\mathrm{vel} (n) 
		\end{equation*}
		for every $m,n \geq 2$.
	\end{proposition}

\begin{proof}[Proof of Proposition~\ref{prop:improved_submult_generalized}]
Let $G^m_{[0,n]}$, $\Gamma^m_\mathrm{vel}(n)$, and
$X_i = \max_{1\leq \ell \leq m} D^\mathrm{vel}_{[im,im+\ell]}(im,im+\ell-1)$
 be as in the proof of Proposition~\ref{prop:improved submult} except that we replace the chemical distance with the vertex extremal length. Since the $X_i$ are independent copies of $\Gamma_\mathrm{vel}(m)$ and are exponentially unlikely to exceed their mean by a large factor by Lemma~\ref{lem:NBU_generalized}, it suffices as in the proof of Proposition~\ref{prop:improved submult} to prove that
\[
\Gamma_\mathrm{vel}(mn) \leq 4 \sum_{i=1}^{\min\{\lceil \Gamma^m_\mathrm{vel}(n)\rceil+1,n\}} X_{(i)},
\]
the rest of the argument being similar. To this end, we fix $p=\ell m+r$ with $0\leq \ell \leq n-1$ and $1\leq r\leq m$, the case $\ell=0$ being trivial, and let $\mathscr{P}$ be a random simple path from $0$ to $\ell-1$ in $G^m_{[0,\ell]}$ attaining the minimum in the definition of the VEL, so that $\sum_i \p(i\in \mathscr{P}\mid \cX)^2 \leq \Gamma^m_\mathrm{vel}(n)$. For each $x,y\in \{im,\ldots,(i+1)m-1\}$, concatenating independent optimal random paths from $im$ to $x$ in $G_{[im,x+1]}$ and from $im$ to $y$ in $G_{[im,y+1]}$ yields a random path from $x$ to $y$ in $G_{[im,(i+1)m]}$ whose hitting probabilities have sum of squares at most $4X_i$. As in the proof of Proposition~\ref{prop:improved submult}, we can therefore turn $\mathscr{P}$ into a random path $\overline{\mathscr{P}}$ from $0$ to $p-1$ in $G_{[0,p]}$ by replacing the visit of $\mathscr{P}$ to each block $i\leq \ell-1$ by such a random path between the relevant entry and exit points (taken conditionally independent given $\mathscr{P}$), and appending the edge $\{\ell m-1,\ell m\}$ and an optimal random path from $\ell m$ to $p-1$ in $G_{[\ell m,p]}$. We now bound the hitting probabilities of $\overline{\mathscr{P}}$. Fix $0\leq i\leq \ell-1$ and a vertex $w$ in the $i$th block $\{im,\ldots,(i+1)m-1\}$. Since $\mathscr{P}$ visits each block at most once, $\overline{\mathscr{P}}$ can only hit $w$ if $i\in\mathscr{P}$, and in this case only during the random path inserted for the visit of $\mathscr{P}$ to block $i$. Let $h(w)$ denote the conditional probability given $\cX$ and $\mathscr{P}$ that this inserted path hits $w$, so that $\sum_{w} h(w)^2\leq 4X_i$ on the event $\{i\in \mathscr{P}\}$, where the sum is over the vertices of block $i$. We have by Cauchy--Schwarz that
\[
\p(w\in\overline{\mathscr{P}}\mid\cX)^2=\E\bigl[\mathbbm{1}(i\in\mathscr{P})\,h(w)\bigm|\cX\bigr]^2\leq \p(i\in\mathscr{P}\mid\cX)\;\E\bigl[\mathbbm{1}(i\in\mathscr{P})\,h(w)^2\bigm|\cX\bigr],
\]
and summing over the vertices $w$ of block $i$ gives that
\[
\sum_{w=im}^{(i+1)m-1} \p(w\in\overline{\mathscr{P}}\mid\cX)^2\leq \p(i\in\mathscr{P}\mid\cX)\;\E\bigl[\mathbbm{1}(i\in\mathscr{P})\, 4X_i\bigm|\cX\bigr] = 4X_i\,\p(i\in\mathscr{P}\mid\cX)^2.
\]
Since the vertices $w\geq \ell m$ can be visited only by the final path from $\ell m$ to $p-1$, they contribute at most $X_\ell$, and it follows that
\[
	D^\mathrm{vel}_{[0,p]}(0,p-1) \leq \sum_{w} \p(w\in \overline{\mathscr{P}} \mid \cX)^2 \leq 4 \sum_{i=0}^{\ell-1} X_i \p(i\in \mathscr{P}\mid \cX)^2 + X_\ell.
\]
Since $\p(i\in \mathscr{P}\mid \cX)^2\in[0,1]$ for every $i$ and $\sum_i \p(i\in \mathscr{P}\mid \cX)^2 \leq \Gamma^m_\mathrm{vel}(n)$,  it follows by the rearrangement inequality that 
\[D^\mathrm{vel}_{[0,p]}(0,p-1)\leq 4\sum_{i=1}^{\min\{\lceil \Gamma^m_\mathrm{vel}(n)\rceil+1,n\}} X_{(i)},\] which implies the claim since $1\leq p \leq mn$ was arbitrary.
\end{proof}

\subsection{Point-to-point resistances}
\label{subsec:point_to_point}

We now adapt the arguments of Section~\ref{subsec:the_diameter} to the  VEL and deduce our desired consequences for the effective resistance. Since the proofs here are very similar to those of Section~\ref{subsec:the_diameter}, we keep the details light and focus on what needs to be changed; the main change is that we work entirely with block graphs instead of passing to the continuum quantities (which we have not defined).
Our first goal is to prove the following analogue of Lemma~\ref{lem:domination_time}. (This lemma is much weaker than Lemma~\ref{lem:domination_time} when $s\ll t$, but this will not cause any problems.)

\begin{lemma}\label{lem:domination_vel}
If $\cX$ satisfies \emph{(S)} then $D^\mathrm{vel}_{[0,s]}(0,s-1)$ is stochastically dominated by $3(t/s+3)D^\mathrm{vel}_{[0,t]}(0,t-1)$ for all integers $1\leq s\leq t$.
\end{lemma}

To prove this lemma we first study how the VEL changes under coarse graining. 
 For each pair of integers $1\leq s\leq t$, let $G^{t/s}_{[0,s]}$ be defined as in the proof of Proposition~\ref{prop:improved submult}, so that $i,j\in\{0,\ldots,s-1\}$ are adjacent if $|i-j|=1$ or there exists $\{x,y\}\in\cX$ with $x\in[(t/s) i,(t/s)(i+1)]$ and $y\in[(t/s) j,(t/s)(j+1)]$.

\begin{lemma}\label{lem:coarse_graining}
For each $\cX\subseteq\Omega$ we have that
\[ \mathscr{V}\Bigl(0\leftrightarrow s-1;G^{t/s}_{[0,s]}\Bigr) \leq 3(t/s+3)\, \mathscr{V}\Bigl(0\leftrightarrow t-1;G_{[0,t]}\Bigr)\]
for every pair of integers $1\leq s\leq t$.
\end{lemma}

\begin{proof}[Proof of Lemma~\ref{lem:coarse_graining}]
Let $\lambda=t/s$.
For each $k\in\{0,\ldots,t-1\}$ let $A(k)$ be the set of $i\in\{0,\ldots,s-1\}$ with $[k,k+1]\cap[\lambda i,\lambda(i+1)]\neq\emptyset$, so that each $A(k)$ consists of at most three consecutive integers and each $i$ belongs to at most $\lambda+3$ of the sets $A(k)$. Observe moreover that $0\in A(0)$, $s-1\in A(t-1)$, and that if $k$ and $k'$ are adjacent in $G_{[0,t]}$ then there exists $i\in A(k)$ and $i'\in A(k')$ that are either equal or adjacent in $G^\lambda_{[0,s]}$. Thus, any path $\mathscr{P}$ from $0$ to $t-1$ in $G_{[0,t]}$ can be converted into a path from $0$ to $s-1$ in $G^\lambda_{[0,s]}$ visiting only vertices in $\bigcup_{k\in\mathscr{P}}A(k)$. Applying this to the random path $\mathscr{P}$ attaining the VEL in $G_{[0,t]}$ and using Cauchy--Schwarz yields that the VEL in $G^\lambda_{[0,s]}$ is at most
\[\sum_i \Biggl[\sum_{k\in A^{-1}(i)}\p(k\in\mathscr{P})\Biggr]^2 \leq (\lambda+3)\sum_i\sum_{k\in A^{-1}(i)}\p(k\in\mathscr{P})^2\leq 3(\lambda+3)\sum_k\p(k\in\mathscr{P})^2\]
as claimed.
\end{proof}

\begin{proof}[Proof of Lemma~\ref{lem:domination_vel}]
This follows immediately from Lemma~\ref{lem:coarse_graining} since $\frac{s}{t}\cX$ has the same law as $\cX$ when $\cX$ satisfies (S).
\end{proof}

We now state and prove the VEL analogue of Lemmas~\ref{lem:regular_times} and \ref{lem:gamma_to_diam}.
We say that $n\geq 1$ is \textbf{$\eta$-regular for the VEL} if $n^{-\eta}\E\Gamma_\mathrm{vel}(n)\geq \eta\, m^{-\eta}\E\Gamma_\mathrm{vel}(m)$ for every $m\leq n$.

\begin{lemma}\label{lem:gamma_to_diam_vel}
Let $\cX$ satisfy \emph{(I)}, \emph{(T)}, and \emph{(S)}. Then Lemma~\ref{lem:regular_times} holds with $\Gamma_\mathrm{vel}$ and $\theta_\mathrm{vel}$ in place of $\gamma$ and $\theta$, and for each $\eta>0$ there exists $c(\eta)>0$ such that
\[\p\bigl(D^\mathrm{vel}_{[0,n]}(0,n-1)\geq c(\eta)\E\Gamma_\mathrm{vel}(n)\bigr)\geq c(\eta)\]
for every $n\geq 1$ that is $\eta$-regular for the VEL.
\end{lemma}

\begin{proof}[Proof of Lemma~\ref{lem:gamma_to_diam_vel}]
The first claim is proven exactly as Lemma~\ref{lem:regular_times}. For the second, the proof of Lemma~\ref{lem:gamma_to_diam} applies with $f(m)=D^\mathrm{vel}_{[0,m]}(0,m-1)$ in place of $d_{[0,z]}(0,z)$, a uniform random element of $\{0,\ldots,n-1\}$ in place of a uniform random point of $[0,t]$, Lemma~\ref{lem:NBU_generalized} in place of the NBU property of $\gamma$, and Lemma~\ref{lem:domination_vel} in place of Lemma~\ref{lem:domination_time}: the only properties of $f$ used there are that $f(m+1)\leq f(m)+1$, that $f(m)\leq f(k)+D^\mathrm{vel}_{[k,m]}(k,m-1)$ for $k<m$, and translation invariance. The factor of $t/s$ in Lemma~\ref{lem:domination_vel} is harmless since the desired lower bound on $D^\mathrm{vel}_{[0,X]}(0,X-1)$ is proven to hold with good probability for a uniform random $X \in \{0,\ldots,n-1\}$, and such a point has $n/X$ of constant order with good probability.
\end{proof}

\begin{proof}[Proof of Theorem~\ref{thm:main_discrete_other_metrics}]
As in the proof of Corollary~\ref{cor:first_exponent}, it follows from Lemma~\ref{lem:NBU_generalized} and Proposition~\ref{prop:improved_submult_generalized} that $\E\Gamma_\mathrm{vel}(n)=n^{\theta_\mathrm{vel}+o(1)}$ for some $\theta_\mathrm{vel}\in [0,1]$. Let $\Delta_n$ be the maximal degree of $G_{[0,n]}$, and let $\eps>0$. Since effective resistances define a metric and are decreasing in the graph, Lemma~\ref{lem:EEL_VEL} yields that
\[ \sup_{x,y} R_{[0,n]}(x,y) \leq 2\Gamma_\mathrm{res}(n) \leq 2\Gamma_\mathrm{vel}(n) \qquad \text{and} \qquad R_{[0,n]}(0,n-1)\geq \frac{1}{2\Delta_n}D^\mathrm{vel}_{[0,n]}(0,n-1).\]
As such, it suffices to prove that $\E[D^\mathrm{vel}_{[0,n]}(0,n-1)\mathbbm{1}(\Delta_n\leq n^\eps)]\geq n^{\theta_\mathrm{vel}-2\eps-o(1)}$, the claim then following with $\theta_\mathrm{res}=\theta_\mathrm{vel}$. This is trivial if $\theta_\mathrm{vel}=0$. Otherwise, fix $0<\eta<\theta_\mathrm{vel}$ and let $m\in[n^{1-\eps},n]$ be $\eta$-regular for the VEL, which exists for large $n$ by Lemma~\ref{lem:gamma_to_diam_vel}. Lemmas~\ref{lem:domination_vel} and \ref{lem:gamma_to_diam_vel} imply that 
\[D^\mathrm{vel}_{[0,n]}(0,n-1)\geq \frac{c(\eta)\,m}{3(n+3m)}\E\Gamma_\mathrm{vel}(m)\geq n^{\theta_\mathrm{vel}-2\eps-o(1)}\] with probability at least $c(\eta)$, and the claim follows since $\p(\Delta_n\leq n^\eps)\to 1$ by (L).
\end{proof}

\subsection{Hypothesis (L) for the Brownian trace}
\label{subsec:L_brownian}
 
In this section we prove that the trace of Brownian motion satisfies hypothesis (L) in two and three dimensions.
 
\begin{lemma}\label{lem:brownian_L}
Let $d\in \{2,3\}$ and let $(B_t)_{t\in\R}$ be a two-sided Brownian motion in $\R^d$ with $B_0=0$. Let $\cX=\{\{s,t\}:B_s=B_t\}$, and let $N=\#\{i\in\Z: B[0,1]\cap B[i,i+1]\neq\emptyset\}$.
\begin{enumerate}
\item If $d=3$ then $\E[N^q]<\infty$ for every $q\geq 1$.
\item If $d=2$ then for each $\eps>0$ there exists a constant $C$ such that the maximal degree of $G_{[0,n]}$ is at most $C(\log n)^3$ with probability at least $1-\eps$ for every $n\geq 2$.
\end{enumerate}
In particular, $\cX$ satisfies \emph{(L)} in both cases.
\end{lemma}
 
\begin{proof}[Proof of Lemma~\ref{lem:brownian_L}]
Letting $N_+(n)= \#\{i\in\{1,\ldots,n\}: B[0,1]\cap B[i,i+1]\neq\emptyset\}$, it suffices to prove that there exists a positive constant $c$ such that
\[
	\p\left(N_+(n) \geq k (1+ \mathbbm{1}(d=2)\log n)\right) \leq e^{-ck^{1/2}}
\]
for every $n,k \geq 2$, with the claims concerning the maximum degree in $G_{[0,n]}$ following by a union bound (applying the same estimate to the time-reversed process $(B_{-t})_{t\geq 0}$ to control the intervals $[i,i+1]$ with $i<0$).

We prove a bound on the tail of the number of integers $0\leq i \leq n$ such that $B[i,i+1]=\{B(t): t\in [i,i+1]\}$ intersects the ball $U_r:=\{x\in \R^d : \|x\|<r\}$, which we denote by $N_+(r,n)$, then use this to estimate the tail of the degree. 
Define two sequences of stopping times $(\sigma_i)_{i\geq 0}$ and $(\tau_i)_{i\geq 0}$ recursively by $\sigma_0=0$, $\tau_i=\inf\{t\geq \sigma_i : |B_t| \geq 2r\}$, and $\sigma_{i+1}=\inf \{t\geq \tau_i: |B_t|\leq r\}$.
 Every time $t\geq 0$ with $B_t\in U_r$ belongs to an interval of the form $[\sigma_i,\tau_i]$ for some $i\geq 0$, while each such interval $[\sigma_i,\tau_i]$ intersects at most $\tau_i-\sigma_i+2$ of the intervals $[j,j+1]$. Thus, letting $K(r,n)=\sup\{k\geq 0:\sigma_k\leq n\}$ for each integer $n\geq 1$, we have that
\begin{equation}
\label{eq:excursion_count}
N_+(r,n)\leq \sum_{k=0}^{K(r,n)}(\tau_k-\sigma_k+2).
\end{equation}
Let $\cF_t$ denote the $\sigma$-algebra generated by $(B_s)_{0\leq s \leq t}$. If $d=3$ then there exists a positive constant $c$ (independent of $r$) such that $\p(\sigma_{i+1}=\infty \mid \cF_{\tau_i})=c>0$ whenever $\tau_i<\infty$ for every $i\geq 0$, while when $d=2$ we have instead that
\[
	\p(\sigma_{i+1}-\tau_i \geq t r^2 \mid \cF_{\tau_i}) \geq \frac{c}{\log t}
\]
for every $t\geq 2$ (this follows by first using the function $\log r^{-2} \|x\|_2^2$, which is harmonic on $\R^2\setminus \{0\}$, to prove that Brownian motion started at distance $2r$ reaches distance $r t$ before hitting $U_r$ with probability of order $1/\log t$, then noting that the Brownian motion is unlikely to return to $U_r$ in time $tr^2$ after reaching distance $rt$). Using the strong Markov property, this easily implies that there exists a positive constant $c_1$ such that $\p( K(r,n) \geq \lambda ) \leq e^{-c_1\lambda} $ when $d=3$ and 
\[\p\left( K(r,n) \geq \lambda \left(1 \vee \log \frac{n}{r^2}\right)\right) \leq e^{-c_1\lambda}\] when $d=2$.
Since a Brownian motion started anywhere in $U_{2r}$ has at least a constant probability to leave $U_{2r}$ in time $r^2$, we also have by the strong Markov property that there exists a positive constant $c_2$ such that
\[
	\p(\tau_i-\sigma_i \geq m r^2 \mid \cF_{\sigma_i}) \leq e^{-c_2 m}
\]
for every integer $m\geq 1$. Letting $m,r,n\geq 2$ and writing 
\[M= \left\lceil m \left(1 \vee \log \frac{n}{r^2} \mathbbm{1}(d=2)\right)\right\rceil,\] it follows by a union bound and standard large deviation estimates that there exist positive constants $c_4$ and $C$ such that
\begin{multline*}
	\p\left(N_+(r,n) \geq Cr^2 M \right)
	\leq \p\left(K(r,n) \geq M\right)
	+ 
	\p\left(\sum_{i=0}^{M-1} (\tau_i -\sigma_i +2) \geq C r^2 M\right)
	\leq e^{-c_4 m}.
\end{multline*}
Now, we also have that by standard properties of Brownian motion that there exists a positive constant $c_5$ such that $\p(\max_{0\leq t \leq 1}\|B_t\| \geq r) \leq e^{-c_5 r^2}$ for every $r\geq 1$. Thus, we have by a union bound that there exists a positive constant $c_6$ such that if $r\leq n^{3/8}$ (which ensures that $\log n$ and $\log n/r^2$ are of the same order) then
\begin{multline*}
	\p(N_+(n) \geq k (1 + \mathbbm{1}(d=2) \log n)) \leq \p(N_+(r,n) \geq k (1 + \mathbbm{1}(d=2) \log n))
	\\ + \p(\max_{0\leq t \leq 1}\|B_t\| \geq r)
	\leq e^{-c_6 k/r^2} + e^{-c_5 r^2},
\end{multline*}
and the claim follows by taking $r=k^{1/4}$.
\end{proof}

\section*{Acknowledgements
}
	TH was supported by NSF grant DMS-1928930 and a Packard Fellowship for Science and Engineering. The use of the duality formula for vertex extremal length was inspired by conversations between TH and Asaf Nachmias circa 2015.
	
\section*{AI usage statement}

AI assistants were used for manuscript review, proofreading, and the creation of figures.

\footnotesize{
\bibliographystyle{abbrv}
\bibliography{unimodularthesis}
}

\end{document}